\documentclass[reqno, oneside]{amsart}
\usepackage{style}
\usepackage{hyperref}
\numberwithin{equation}{section}

\begin{document}
\title[Szeg\H{o} theorem]{A spectral Szeg\H{o} theorem on the real line}
\author{R.~V.~Bessonov, \quad  S.~A.~Denisov}

\address{
\begin{flushleft}
Roman Bessonov: bessonov@pdmi.ras.ru\\\vspace{0.1cm}
St.\,Petersburg State University\\  
Universitetskaya nab. 7/9, 199034 St.\,Petersburg, RUSSIA\\
\vspace{0.1cm}
St.\,Petersburg Department of Steklov Mathematical Institute\\ Russian Academy of Sciences\\
Fontanka 27, 191023 St.Petersburg,  RUSSIA\\
\end{flushleft}
}
\address{
\begin{flushleft}
Sergey Denisov: denissov@wisc.edu\\\vspace{0.1cm}
University of Wisconsin--Madison\\  Department of Mathematics\\
480 Lincoln Dr., Madison, WI, 53706, USA\vspace{0.1cm}\\
\vspace{0.1cm}
Keldysh Institute of Applied Mathematics\\ Russian Academy of Sciences\\
Miusskaya pl. 4, 125047 Moscow, RUSSIA\\
\end{flushleft}
}

\thanks{The first author is supported by RFBR grant mol\_a\_dk 16-31-60053. The work of the second author on the first two sections was supported by RSF-14-21-00025 and his research on the rest of the paper was supported by NSF Grant DMS-1464479.}
\subjclass[2010]{42C05, 34L40, 34A55}
\keywords{Szeg\H{o} class, Canonical Hamiltonian system, String equation,  Inverse problem, Entropy, Muckenhoupt class}

\begin{abstract}
We characterize even measures $\mu = w\, dx + \mu_s$ on the real line $\R$ with finite entropy integral  $\int_{\R}\frac{\log w(t)}{1+t^2}\,dt > -\infty$ in terms of $2 \times 2$ Hamiltonians generated by $\mu$ in the sense of the inverse spectral theory. As a corollary, we obtain  criterion for spectral measure of Krein string to have converging logarithmic integral.
\end{abstract}

\maketitle
\large

\section{Introduction}\label{s1}
Each probability measure $\mu$ supported on an infinite subset of the unit circle $\T = \{z: \, |z| = 1\}$ of the complex plane, $\C$, gives rise to the infinite family $\{\Phi_n\}_{n \ge 0}$ of monic polynomials orthogonal with respect to~$\mu$. For integer $n \ge 0$, the polynomial $\Phi_n$ has degree $n$, unit coefficient in front of~$z^n$, and $(\Phi_n,\Phi_k)_{L^2(\mu)} = 0$ for all $k\neq n$. The polynomials $\{\Phi_n\}_{n \ge 0}$ satisfy the recurrence relation 
\begin{equation}\label{eq91}
\Phi_{n+1}(z) = z \Phi_n(z) -\bar\alpha_n \Phi_n^\ast(z), \qquad \Phi_0 =1,
\end{equation}
where $\{\Phi_n^\ast\}$ are the ``reversed'' polynomials defined by $\Phi_n^\ast(z) = z^n\ov{\Phi_n(1/\bar z)}$. Recurrence coefficients $\{\alpha_n\}$ are completely determined by $\mu$ and we have $|\alpha_n|<1$ for every $n \ge 0$. Given any sequence of complex numbers $\{\alpha_n\}$ with $|\alpha_n|<1$, one can find the unique probability measure $\mu$ on $\T$ such that $\{\alpha_n\}$ is the sequence of the recurrence coefficients of $\mu$, see  \cite{Simonbook}, \cite{Szego}.   
\begin{SzegoTheorem}
Let $\mu = w\,dm + \mu_s$ be a probability measure on $\T$ with density $w$ and a singular part $\mu_s$ with respect to the Lebesgue measure $m$ on $\T$. The following assertions are equivalent:
\begin{itemize}
\item[$(a)$] the set $\spn\{z^n,\; n \ge 0\}$ of analytic polynomials is not dense in $L^2(\mu)$;
\item[$(b)$] the entropy of $\mu$ is finite: $\int_{\T}\log w\,dm > -\infty$;
\item[$(c)$] the recurrence coefficients $\{\alpha_n\}$ of $\mu$ satisfy $\sum_{n\ge 0} |\alpha_n|^2 < \infty$.
\end{itemize}
\end{SzegoTheorem}
We refer the reader to \cite{Simonbook}, \cite{SimonDes} for the historical account and an extended version of this result. Independent contributions to different aspects of its proof were done by Szeg\H{o}, Verblunsky, and Kolmogorov. A partial counterpart  of  Szeg\H{o} theorem for measures supported on the real line,~$\R$, is due to Krein \cite{Krein45} and Wiener \cite{Wiener49} (see also Section 4.2 in \cite{DymMcKean} or Theorem A.6 in \cite{Den06} for modern expositions). Denote by $\Pi(\R)$ the class of all Radon measures on $\R$ such that $\int_{\R}\frac{d\mu(t)}{1+t^2} < \infty$. 
\begin{KWtheorem} Let $\mu = w\,dx + \mu_s$ be a measure in $\Pi(\R)$ where $w$ is the density with respect to the Lebesgue measure $dx$ on $\R$ and  $\mu_s$ is the singular part. The following assertions are equivalent:
\begin{itemize}
\item[$(a)$] the set of functions whose Fourier transform is smooth and compactly supported on $[0, +\infty)$ is not dense in $L^2(\mu)$;
\item[$(b)$] the entropy of $\mu$ is finite: $\int_{\T}\frac{\log w(t)}{1+t^2}\,dt > -\infty$.
\end{itemize}
\end{KWtheorem}
Szeg\H{o} and Krein-Wiener theorems have  probabilistic interpretation. Roughly, it says that a stationary Gaussian sequence/process with the spectral measure $\mu$ is non-deterministic if and only if the entropy of $\mu$ is finite, see, e.g, Section II.2 in \cite{ibr} or survey \cite{Bingham} for more details.

\medskip

The aim of this paper is to complement assertions $(a)$, $(b)$ in Krein--Wiener theorem with a necessary and sufficient condition similar to condition $(c)$ in Szeg\H{o} theorem. Instead of recurrence relation $\Phi_{n+1}(z) = z \Phi_n(z) -\bar\alpha_n \Phi_n^\ast(z)$, we will consider canonical Hamiltonian system $JM' = z\Hh M$ which naturally appears from $\mu$ via Krein\,--\,de\,Branges spectral theory.   

\medskip

\noindent Consider the Cauchy problem for a canonical Hamiltonian system on the half-axis $\R_+ = [0, +\infty)$,
\begin{equation}\label{eq1}
J M'(t,z) = z\Hh(t) M(t,z), \qquad M(0,z) = \idm, \quad t\ge 0, \quad z\in \C.
\end{equation}
Here $J = \jm$, the derivative of $M$ is taken with respect to $t$, the Hamiltonian $\Hh$ is the mapping taking numbers $t \in \R_+$ into positive semi-definite matrices, the entries of $\Hh$ are real measurable functions on $\R_+$ absolutely integrable on compact subsets of $\R_+$. In addition, we assume that the trace of $\Hh$ does not vanish identically on any set of positive Lebesgue measure. The Hamiltonian $\Hh$ is called singular if $\int_{0}^{\infty}\trace\Hh(t)\,dt = +\infty$. We say that $\Hh$ is nontrivial if there is no subset $E \subset \R_+$ of full Lebesgue measure such that $\Hh = \left(\begin{smallmatrix}1&0\\0&0\end{smallmatrix}\right)$ on $E$ or $\Hh = \left(\begin{smallmatrix}0&0\\0&1\end{smallmatrix}\right)$ on $E$.

\medskip

Let $\Hh$ be a singular nontrivial Hamiltonian on $\R_+$, and let $M$ be the solution of \eqref{eq1}. Fix a parameter $\omega \in \R \cup \{\infty\}$ and define the Weyl-Titchmarsh function $m$ of \eqref{eq1} on $\C\setminus\R$ by
\begin{equation}\label{eq9}
m(z) = \lim_{t \to +\infty}\frac{\omega\Phi^+(t,z) + \Phi^-(t,z)}{\omega\Theta^+(t,z) + \Theta^-(t,z)}, \qquad
M(t, z) =  \Bigl(\!\begin{smallmatrix}\Theta^+(t,z) & \Phi^+(t,z) \\ \Theta^-(t,z) & \Phi^-(t,z) \end{smallmatrix}\!\Bigr). 
\end{equation}
The fraction $\frac{\infty c_1 + c_2}{\infty c_3 + c_4}$ for non-zero numbers $c_1$, $c_3$ is interpreted as $\frac{c_1}{c_3}$. 
For the Weyl-Titchmarsh theory of canonical Hamiltonian systems see \cite{HSW} or Section 8 in \cite{Romanov}. Theorem~2.1 in \cite{HSW} implies that the denominator of the fraction in \eqref{eq9} is nonzero for large $t \ge 0$, the function $m$ does not depend on the choice of the parameter $\omega$, and $\Im m(z) > 0$ for $z$ in $\C^+ = \{z \in \C: \; \Im z >0\}$. Hence, there exists a measure $\mu \in \Pi(\R)$, and numbers $a \in \R$, $b \ge 0$, such that 
\begin{equation}\label{eq371}
m(z) = \frac{1}{\pi}\int_{\R}\left(\frac{1}{x-z} - \frac{x}{1+x^2}\right)d\mu(x) + bz + a, \qquad z \in \C \setminus \R. 
\end{equation}
The measure $\mu$ in \eqref{eq371} is called the spectral measure of system \eqref{eq1}. Two singular nontrivial Hamiltonians $\Hh_1$, $\Hh_2$ on~$\R_+$ are called equivalent if there exists an increasing absolutely continuous function $\eta$ on $\R_+$, $\eta(0) = 0, \lim_{t\to\infty}\eta(t)=\infty$, such that $\Hh_2(t) = \eta'(t)\Hh_1(\eta(t))$ for Lebesgue almost all $t \in \R_+$. It is easy to check that equivalent Hamiltonians have equal Weyl-Titchmarsh functions, see \cite{WW12}. 
The following theorem is central to Krein -- de Branges inverse spectral theory \cite{KK68}, \cite{dbbook}.
\begin{DBtheorem} For every analytic function $m$ in $\C^+$ with positive imaginary part, there exists a singular nontrivial Hamiltonian $\Hh$ on $\R_+$ such that $m$ is the Weyl-Titchmarsh function~\eqref{eq9} for $\Hh$. Moreover, any two singular nontrivial Hamiltonians $\Hh_1$, $\Hh_2$ on~$\R_+$ generated by $m$ are equivalent. 
\end{DBtheorem}
See \cite{Romanov}, \cite{Winkler95} for proofs to this theorem. 
%
%
A measure $\mu$ on $\R$ is called even if $\mu(I) = \mu(-I)$ for every interval $I \subset \R_+$. It is well-known that a Hamiltonian $\Hh$ has the diagonal form $\Hh = \diag(h_1, h_2)$ almost everywhere on $\R_+$ if and only if its spectral measure $\mu$ is even and $a = 0$ in \eqref{eq371}, see Lemma~\ref{l11} below. Here $\diag(c_1, c_2) = \left(\begin{smallmatrix}c_1 & 0 \\ 0 &c_2 \end{smallmatrix}\right)$ for $c_1$, $c_2 \in \R_+$.

\medskip

Szeg\H{o} class $\sz$ on the real line $\R$ consists of measures $\mu \in \Pi(\R)$ that satisfy equivalent assertions $(a)$, $(b)$ in Krein--Wiener theorem. Given a measure $\mu = w\,dx + \mu_s$ in $\sz$, define its normalized entropy by   
$$
\K(\mu) = \log \frac{1}{\pi}\int_{\R}\frac{d\mu(x)}{1+x^2} - \frac{1}{\pi}\int_{\R}\frac{\log w(x)}{1+x^2}\,dx.
$$
By Jensen inequality, we have $\K(\mu) \ge 0$, and, moreover, $\K(\mu) = 0$ if and only if $\mu$ is a non-zero scalar multiple of the Lebesgue measure on $\R$. 

\medskip

We say that a  measure $\mu \in \Pi(\R)$ generates a Hamiltonian $\Hh$ if the Weyl-Titchmarsh function~\eqref{eq9} of $\Hh$ has the form $m: z \mapsto \frac{1}{\pi}\int_{\R}\bigl(\frac{1}{x-z} - \frac{x}{1+x^2}\bigr)\,d\mu(x)$. To every $\Hh$ with $\sqrt{\det\Hh} \notin L^1(\R_+)$ we associate the sequence of points $\{\eta_n\}$ by
\begin{equation}\label{eq71}
\eta_n = \min\left\{t \ge 0: \int_{0}^{t}\sqrt{\det\Hh(s)}\,ds = n\right\}, \qquad n \ge 0.
\end{equation}
Our main result is the following theorem.
\begin{Thm}\label{t1}
An even measure $\mu \in \Pi(\R)$ belongs to the Szeg\H{o} class $\sz$ if and only if some (and then every)  Hamiltonian $\Hh = \diag(h_1, h_2)$ generated by $\mu$ is such that $\sqrt{\det \Hh} \notin L^1(\R_+)$ and 
\begin{equation}\label{eq2}
\widetilde \K(\Hh) = \sum_{n = 0}^{+\infty}  \left(\int_{\eta_n}^{\eta_{n+2}}h_1(s)\,ds \cdot \int_{\eta_n}^{\eta_{n+2}}h_2(s)\,ds - 4\right) < \infty,
\end{equation}
where $\{\eta_n\}$ are given by \eqref{eq71}. 
Moreover, we have $\widetilde\K(\Hh) \le c\K(\mu)e^{c\K(\mu)}$ and $\K(\mu) \le c\widetilde\K(\Hh)e^{c\widetilde\K(\Hh)}$ for an absolute constant $c$. 
\end{Thm}
\noindent By definition, the terms in \eqref{eq2} are nonnegative:
$$
\int_{\eta_{n}}^{\eta_{n+2}}h_1(s)\,ds\cdot \int_{\eta_{n}}^{\eta_{n+2}}h_2(s)\,ds - 4 \ge \left(\int_{\eta_n}^{\eta_{n+2}}\sqrt{\det\Hh(s)}\,ds\right)^2 - 4 = 0,  
$$
and the sum in \eqref{eq2} equals zero if and only if $\Hh$ is a constant Hamiltonian. Note that the spectral measure $\mu$ of a constant diagonal Hamiltonian $\Hh$ with $\det\Hh \neq 0$ is a scalar multiple of the Lebesgue measure on $\R$, in particular, we have $\K(\mu) = 0$ in this case.

\medskip

Diagonal canonical Hamiltonian systems are closely related to the differential equation of a vibrating string:
\begin{equation}\label{eq40}
-\frac{d}{dM(t)}\frac{d}{dt}\Bigl( y(t, z)\Bigr) = z y(t, z), \qquad t \in [0, L), \qquad z \in \C. 
\end{equation} 
Here $0 < L \le +\infty$ is the length of the string, $M: (-\infty, L) \to \R_+$ is an arbitrary non-decreasing and right-continuous function (mass distribution) that satisfies $M(t)=0$ for $t<0$. If $M$ is smooth and strictly increasing on $\mathbb{R}_+$, then equation \eqref{eq40} takes the form $-y'' = z M' y$. 

\medskip

\noindent In this paper, we consider $L$ and $M$ that satisfy the following conditions:
\begin{equation}\label{sdsing}
L+\lim_{t\to L}M(t)=\infty\,\quad {\rm and}\,\quad \lim_{t\to L}M(t)>0\,,
\end{equation}
where the last bound means that $M$ is not identically equal to zero. If \eqref{sdsing} holds, we say that $M$ and $L$ form $[M,L]$ pair.
 To every  $[M, L]$ pair one can associate a string and Weyl-Titchmarsh function $q$ with spectral measure $\sigma$ supported on the positive half-axis $\R_+$. We discuss these objects in more detail in Section~\ref{s6}.  Theorem \ref{t1} can be reformulated for Krein strings  as follows.
\begin{Thm}\label{t2}
Let $[M, L]$ satisfy \eqref{sdsing} and $\sigma = v\,dx + \sigma_s$ be the spectral measure of the corresponding string. Then, we have
$\int_{0}^{\infty}\frac{\log v(x)}{(1+x)\sqrt{x}}\,dx >- \infty$ 
if and only if $\sqrt{M'} \notin L^1(\R_+)$ and 
\begin{equation}\label{sdk7}
\widetilde\K[M,L] = \sum_{n=0}^{+\infty}\Bigl((t_{n+2} - t_n)(M(t_{n+2}) - M(t_n)) - 4\Bigr) < \infty,
\end{equation}
where $t_n = \min\bigl\{t \ge 0:\; n = \int_{0}^{t}\sqrt{M'(s)}\,ds\bigr\}$.
\end{Thm}
Condition \eqref{sdsing} guarantees that the string $[M,L]$ has the unique spectral measure. It does not restrict the generality of Theorem \ref{t2}: if \eqref{sdsing} is violated, then either $M=0$ and  $\int_{0}^{\infty}\frac{\log v(x)}{(1+x)\sqrt{x}}\,dx = -\infty$ because $v=0$, or $L+\lim_{t\to L} M(t)<\infty$ in which case the Weyl-Titchmarsh function is meromorphic and real-valued on $\mathbb{R}$, so $v(x)=0$ again and the logarithmic integral diverges. 
More details on Theorem \ref{t2} can be found in Section \ref{s6}.

\medskip 

\noindent{\bf Historical remarks.} Except for Krein--Wiener theorem, all previously known results on Szeg\H{o} theorem in the continuous setting were proved for the so-called Krein systems, i.e., differential systems that appear as a result of  ``orthogonalization process with continuous parameter'' invented by Krein in \cite{Kr81}. Krein systems with locally summable coefficients can be reduced to the canonical Hamiltonian systems  with absolutely continuous Hamiltonians $\Hh$ (see, e.g, \cite{B17} for this reduction in the diagonal case). The class of Hamiltonians considered in  Theorem \ref{t1} is considerably wider. Krein himself formulated a restricted version of Szeg\H{o} theorem for Krein systems in \cite{Kr81}. In \cite{Den02}, the second author of this paper characterized  Krein systems with coefficients from a Stummel class whose spectral measures belong to $\sz$. In \cite{Tep05}, Teplyaev  fixed an error in the original formulation of Szeg\H{o} theorem  in \cite{Kr81}. The reader interested in Szeg\H{o} theory for Krein  systems can find further information in monograph \cite{Den06}. In \cite{KS03} and \cite{KS09}, Killip and Simon  proved analogs of Szeg\H{o} theorem for Jacobi matrices and Schr\"odinger operators. See also the work \cite{NPVY} by Nazarov, Peherstorfer, Volberg, and Yuditskii for a closely related subject of sum rules for Jacobi matrices.  

\medskip

\noindent{\bf The structure of the paper.}  We start with studying the basic properties  of entropy function for diagonal canonical systems in Section 2. Section 3 contains the proof of upper and lower bounds for the entropy.  Theorem \ref{t1} is proved in the fourth section. The new functional class which appears in the proof of Theorem \ref{t1} is studied in Section 5. We consider Krein strings and prove Theorem \ref{t2} in Section 6. The paper end with appendix which contains some auxiliary results.\medskip

\noindent{\bf Notation.} In the text, we use the following standard notation. Given set $E\subset \mathbb{R}$ with positive Lebesgue measure $|E|>0$ and nonnegative $f\in L^1(E)$, we denote $\langle f\rangle_E=\frac{1}{|E|}\int_E fdx$. Suppose $a\in \R$, $l>0$, then $I_{a,l}=[a,a+l)$.
The symbols $C,c$ denote absolute constants which can change the value from formula to formula. For two non-negative functions
$f_{1}$, $f_{2}$, we write $f_1\lesssim f_2$ if  there is an absolute constant $C$ such that $f_1\le Cf_2$ for all values of the arguments of 
$f_{1}$, $f_{2}$. We define $\gtrsim$ similarly and say that $f_1\sim f_2$ if $f_1\lesssim f_2$ and $f_2\lesssim f_1$ simultaneously. Given a set $E\subset \mathbb{R}$, $\chi_E$ stands for the characteristic function of $E$. The norm of the space $L^p(\R_+)$ is denoted by $\|\cdot\|_{p}$. The space $L^{1}_{\rm loc}(\R_+)$ consists of functions that are absolutely integrable on compact subsets of $\R_+$.  Symbol $[x]$ stands for the integer part of a real number $x$.\bigskip

\section{Entropy function of a canonical Hamiltonian system}\label{s2}
In this section we introduce the entropy function of a diagonal canonical Hamiltonian system and show that it has a number of remarkable properties. 

\medskip

Let $\Hh = \diag(h_1, h_2)$ be a singular nontrivial diagonal Hamiltonian on $\R_+$, and let $m$, $\mu$ be its Weyl-Titchmarsh function and the spectral measure, so that
\begin{equation}\label{eq372}
\Im m(z) = \frac{1}{\pi}\int_{\R}\frac{\Im z}{|x-z|^2}d\mu(x) + b \Im z, \qquad z \in \C^+. 
\end{equation}

\medskip

For every $r \ge 0$ define $\Hh_r$ to be the Hamiltonian on $\R_+$ taking $x$ into $\Hh(x+r)$. Let $m_r$, $\mu_r$, $b_r$ denote the Weyl-Titchmarsh function, the spectral measure, and the coefficient in \eqref{eq371} of system \eqref{eq1} for $\Hh=\Hh_r$. Each time we work with these objects later in the text we assume that $\Hh_r$ is nontrivial. Define 
\begin{equation}\label{eq63}
\I_{\Hh}(r) = \frac{1}{\pi}\int_{\R}\frac{d\mu_r(x)}{1+x^2} + b_r = -i m_r(i), \qquad
\J_{\Hh}(r) = \frac{1}{\pi}\int_{\R}\frac{\log w_r(x)}{1+x^2}\,dx,
\end{equation}  
where $w_r$ is the density of the absolutely continuous part of $\mu_r = w_r \,dx + \mu_{r,s}$. The second identity above follows from the fact that $\mu$ is even, hence $m$ takes imaginary values on imaginary axis. If $\mu_r \notin \sz$, we put $\J_{\Hh}(r) = -\infty$. Define the entropy function of $\Hh$ by
$$
\K_\Hh(r) = \log \I_{\Hh}(r) - \J_{\Hh}(r), \qquad r \ge 0. 
$$
Note again that Jensen inequality and an estimate $b_r\ge 0$ give
\begin{equation}\label{sdp1}
\K_\Hh(r)\ge 0\,.
\end{equation}
For the ``dual'' Hamiltonian $\Hh^d = J^* \Hh J =  \diag(h_2, h_1)$ we denote the corresponding  objects by $\Hh_{r}^{d}$, $m_{r}^{d}$, $\mu_{r}^{d}$, $b_r^d$, $w_{r}^{d}$, $\I_{\Hh_r^d}$, $\J_{\Hh_r^d}$, and $\K_{\Hh^{d}}$. Note that a Hamiltonian $\Hh$ is singular and nontrivial if and only if $\Hh^d$ is  singular and nontrivial. We also will need the Hamiltonian
\begin{equation}\label{cook59}
\widehat\Hh_r(t) = 
\begin{cases}
\Hh(t), \quad &t \in [0,r),\\
\diag(\I_{\Hh}^{-1}(r), \I_{\Hh}(r)), \quad  &t \in [r, + \infty),
\end{cases}
\end{equation}
which plays the role of ``Bernstein-Szeg\H{o} approximation'' to $\Hh$. From formula \eqref{eq63} we see that the Hamiltonian $\widehat\Hh_r$ is correctly defined  and nontrivial if and only if $m_r(i) \neq 0$, that is, $\Hh_r$ is nontrivial.  Later we will use notation $\widehat\mu_r$ for the spectral measure generated by $\widehat{\Hh_r}$.

\medskip

An analytic function $f$ in the upper half-plane $\C^+ = \{z \in \C: \; \Im z>0\}$ is said to have bounded type if $f = \frac{f_1}{f_2}$ for some bounded analytic functions $f_1$, $f_2$ in $\C^+$, where $f_2$ is not identically zero. Denote by $\N(\C^+)$ the class of all functions of bounded type in $\C^+$. For every function $f \in \N(\C^+)$ we have
\begin{equation}\label{eq12}
\int_{\R}\frac{\bigl|\log |f(x)|\bigr|}{1+x^2}\,dx < \infty,
\end{equation} 
see, e.g., Theorem 9 in \cite{dbbook}. The mean type of a function $f \in \N(\C^+)$ is defined by 
$$
\type_+(f) = \limsup\limits_{y\to+\infty}\frac{\log|f(iy)|}{y}.
$$ 
The upper limit above is finite for every nonzero function $f \in \N(\C^+)$ by Theorem 10 in \cite{dbbook}. A remarkable fact of the spectral theory of canonical Hamiltonian systems is that for every $t \ge 0$ the entries of solution $M(t,z)$ to Cauchy problem \eqref{eq1} are entire functions in $z$ of bounded type in~$\C^{+}$ and their mean type in $\C^+$ equals
\begin{equation}\label{eq7}
\xi_{\Hh}(t) = \int_{0}^{t}\sqrt{\det\Hh(s)}\,ds.
\end{equation} 
This formula has been found by Krein \cite{Krein54} in the setting of the string equation and then proved in full generality by de Branges, see Theorem X in \cite{dbii}. A short proof of~\eqref{eq7} is in Section 6 of \cite{Romanov}.  As a consequence, we have the following result.
\begin{Prop}\label{p2}
Let $\Hh$ be a Hamiltonian on $\R_+$ and let entire function $f(z)$ be one of the entries $\{\Theta^{\pm}(t,z),\Phi^{\pm}(t,z)\}$ of the matrix 
$M$ in \eqref{eq9}. Then, if $f$ in not equal to zero identically, we have 
\begin{equation}\label{eq68}
\frac{1}{\pi}\int_{\R}\log|f(x)|\frac{\Im z}{|x - z|^2}\,dx = \log |f(z)| - \xi_{\Hh}(t)\Im z
\end{equation}
for every $z \in \C_+$.
\end{Prop}
\beginpf 
Take  $f$  as one of $\{\Theta^{\pm}\}$ and denote by $\Theta=\left(\!\begin{smallmatrix}\Theta^+\\ \Theta^-\end{smallmatrix}\!\right)$ the first column of  $M$ in \eqref{eq1}. For every $t>0$, we take inner product of  \eqref{eq1} with $\Theta$ and integrate to get the following well-known identity:
\begin{equation}\label{eq61}
\Im(\Theta^+(t,z)\ov{\Theta^{-}(t,z)}) = \Im z \cdot \int_{0}^{t}\langle \Hh(s) \Theta(s,z), \Theta(s,z) \rangle_{\C^2}\,ds, \qquad z \in \C, 
\end{equation}
where the inner product in $\C^2$ is given by $\langle \left(\begin{smallmatrix}c_1\\ c_2\end{smallmatrix}\right), \left(\begin{smallmatrix}c_3\\ c_4\end{smallmatrix}\right) \rangle_{\C^2} = c_1 \ov{c_3} + c_2 \ov{c_4}$.
This identity implies that either $f$ is identically zero or $f(z) \neq 0$ for $z \in \C\setminus\R$. Function  $f \in \N(\C^+)$, it is smooth on $\R$, and has no zeros in $\mathbb{C}^+$.  So, there exists an outer function $F$ on $\C^+$ such that $f(z) = e^{-i\xi_{\Hh}(r)z}F(z)$, $z \in \C^+$, see Theorem 9 in \cite{dbbook}. Now \eqref{eq68} follows from the mean value theorem for the harmonic function $\log|F|$.  The proof for $\Phi^{\pm}$ is similar. \qed

\medskip

\begin{Prop}\label{p4}
Let $f$ be an analytic function in $\C^+$ such that $\Im f(z)>0$ for all $z \in \C^+$. Then for almost all $x \in \R$ there exists  finite non-tangential limit $f(x) = \lim\limits_{\substack{|z-x|<2 \Im z \\ z \to x }}f(z)$ and 
$$
\frac{1}{\pi}\int_{\R}\log|f(x)|\frac{\Im z}{|x - z|^2}\,dx = \log |f(z)|
$$
for every $z \in \C^+$, where  integral in the left hand side converges absolutely.
\end{Prop}
\beginpf Combine Corollary~4.8 in Section 4 with Exercise 13 in Section 7 of Chapter II in \cite{Garnett}. \qed\bigskip

\medskip

For every $\phi \in [0,\pi)$, set $e_\phi = \left(\begin{smallmatrix}\cos\phi\\\sin\phi\end{smallmatrix}\right)$. An open interval $I \subset \R_+$ is called indivisible for $\Hh$ of type $\phi$ if there is a function $h$ on $I$ such that $\Hh(x)=h(x) e_\phi e_\phi^{\top}$ for almost all $x \in I$, and $I$ is the maximal open interval having this property. Note that a Hamiltonian $\Hh$ on $\R_+$ is nontrivial if $(0, +\infty)$ is not an indivisible interval of type $\phi = 0$ or $\phi = \pi/2$ for $\Hh$. 

\medskip

The following four lemmas are known. We give their proofs in Appendix for the reader's convenience.
\begin{Lem}\label{l5}
Let $\Hh$ be a Hamiltonian on $\R_+$ such that $(0, \ell)$ is indivisible interval of type $\phi \in [0, \pi)$ for $\Hh$. Then the solution $M$ of \eqref{eq1} has the form $M(t,z) = \idm - zJ \int_0^t \Hh(\tau)d\tau$ for every $t \in [0, \ell]$. In particular, for $\Hh = \diag(h_1, h_2)$ and $t \in [0, \ell]$ we have
$$
M(t,z) = 
\begin{cases}
\left(\begin{smallmatrix}1&0\\-z\int_{0}^{t}h_1(s)\,ds&1\end{smallmatrix}\right) &\mbox{if }\phi = 0,\\
\left(\begin{smallmatrix}1&z\int_{0}^{t}h_2(s)\,ds\\0&1\end{smallmatrix}\right) &\mbox{if }\phi = \pi/2
\end{cases}.
$$
\end{Lem} 
\begin{Lem}\label{l11}
Let $\Hh$ be a singular nontrivial Hamiltonian on $\R_+$, and let $m$ be its Weyl-Titchmarsh function \eqref{eq9}. Then, $\Hh$ is diagonal if and only if the measure $\mu$ is even and $a = 0$ in the Herglotz representation \eqref{eq371} of $m$.   
\end{Lem}
\begin{Lem}\label{l12}
Let $\Hh$ be a singular nontrivial Hamiltonian on $\R_+$ and let $m$ be its Weyl-Titchmarsh function. Then, we have $b>0$ in the Herglotz representation \eqref{eq371} of $m$ if and only if $(0,\eps)$ is indivisible interval for $\Hh$ of type $\pi/2$ for some $\eps>0$. Moreover, we have $b = \int_{0}^{\eps}\langle \Hh(t) \zo, \zo\rangle\,dt$ in the latter case. 
\end{Lem}
\begin{Lem}\label{l9}
Let $\Hh = \diag(a_1, a_2)$ be the constant Hamiltonian on $\R_+$ generated by positive numbers $a_1$, $a_2$. Then for all $r \ge 0$ we have $w_r = \sqrt{a_2/a_1}$ on $\R$ and 
\begin{equation}\label{sdp3}
\log\I_{\Hh}(r) = \J_{\Hh}(r) = \log\sqrt{a_2/a_1}\,.  
\end{equation}
\end{Lem}

\medskip

The following lemma is crucial for our paper. 
\begin{Lem}\label{l1}
Let $\Hh = \diag(h_1, h_2)$ be a singular nontrivial Hamiltonian on $\R_+$ and let $\mu$ be the spectral measure of system~\eqref{eq1} generated by $\Hh$. Assume that $\mu \in \sz$. Then for every $r \ge 0$ we have
\begin{itemize}
\item[$(a)$] $\mu_r \in \sz$ and $\mu_r^d \in \sz$,
\item[$(b)$] $\J_\Hh(r)= \J_{\Hh}(0) - 2\xi_{\Hh}(r) + 2\log|\Theta^+(r,i) +i\I_\Hh(r)\Theta^-(r,i)|$,
\item[$(c)$] $\I_{\Hh}(r) = 1/\I_{\Hh^d}(r)$, 
\item[$(d)$] $\K_{\Hh}(r) = \K_{\Hh^d}(r)$,
\item[$(e)$] $\widehat\mu_r \in \sz$ and $\K_{\Hh}(0) = \K_{\widehat\Hh_r}(0) + \K_{\Hh}(r)$, 
\end{itemize}
where $\xi_\Hh$ is defined in \eqref{eq7}.
\end{Lem} 
\beginpf Take $r \ge 0$ and consider solutions 
\begin{equation}\label{eq65}
M(t,z)=\left(\!\begin{smallmatrix}\Theta^+(t,z) & \Phi^+(t,z) \\ \Theta^-(t,z) & \Phi^-(t,z)\end{smallmatrix}\!\right), \qquad 
M_r(t,z) = \left(\!\begin{smallmatrix}\Theta_r^+(t,z) & \Phi_r^+(t,z) \\ \Theta_r^-(t,z) & \Phi_r^-(t,z)\end{smallmatrix}\!\right),
\end{equation} 
of Cauchy problem \eqref{eq1} for the Hamiltonians $\Hh$ and $\Hh_r: x \mapsto \Hh(r+x)$, respectively. We have
\begin{equation}\label{eq3}
M_0(t,z) = M_r(t-r, z)M_0(r,z), \qquad t\ge r, \quad z \in \C.
\end{equation}
Indeed, the right hand side of the above equality satisfies equation $JM' = z \Hh M$ on $[r, \infty)$ and coincides with 
$M_0(t,z)$ at $t = r$.  Multiplying matrices in \eqref{eq3} and using \eqref{eq9} with $\omega=0$, we obtain
\begin{equation}\label{eq62}
m_0(z) = \lim_{t \to +\infty}\frac{\Theta^-_{r}(t-r,z)\Phi^+(r,z) + \Phi^-_{r}(t-r,z)\Phi^-(r,z)}{\Theta^{-}_{r}(t-r,z)\Theta^+(r,z) + \Phi^{-}_{r}(t-r,z)\Theta^-(r,z)},
\end{equation}
Suppose there is $c>0$ such that $(c, +\infty)$ is the indivisible interval of type $\pi/2$ for $\Hh$. Then from Lemma \ref{l5} and formula \eqref{eq62}  we see that $m_0(z) = \frac{\Phi^-(c,z)}{\Theta^-(c,z)}$ for all $z\in \C^+$. Since functions $\Phi^-$, $\Theta^{-}$ are real on the real axis, this implies that $\mu$ is a discrete measure concentrated at zeros of entire function $z \mapsto \Theta^-(c, z)$. In particular, we cannot have $\mu \in \sz$. A similar argument applies in the case where $(c, +\infty)$ is the indivisible interval of type $0$ for some $c >0$. It follows that the Hamiltonian $\Hh_r$  is nontrivial for every $r \ge 0$, in particular, its  Weyl-Titchmarsh function $m_r$ is correctly defined and nonzero. Using \eqref{eq62} and \eqref{eq9}  with $\omega=0$ for $m_r$, we get the relation
\begin{equation}
m_0(z) = \frac{\Phi^+(r,z) + m_r(z)\Phi^-(r,z)}{\Theta^+(r,z) + m_r(z)\Theta^-(r,z)}, \qquad z \in \C^+, \quad r \ge0. \label{eq10}
\end{equation}
Hence,
\begin{align*}
\Im m_0(z) 
= &\frac{\Im\bigl(\Phi^+(r,z)\ov{\Theta^+(r,z)} + |m_{r}(z)|^2 \Phi^-(r,z)\ov{\Theta^-(r,z)}\bigr)}{|\Theta^+(r,z) + m_r(z)\Theta^-(r,z)|^{2}} \\
&+ \frac{\Im\bigl(m_r(z)\bigl(\ov{\Theta^+(r,z)}\Phi^-(r,z) - \Theta^-(r,z)\ov{\Phi^+(r,z)}\bigr) \bigr)}{|\Theta^+(r,z) + m_r(z)\Theta^-(r,z)|^{2}}.
\end{align*}
Since the analytic function $m_r$ has positive imaginary part in $\C^+$  for every $r \ge 0$, we can take non-tangential limit as $z \to x$ in this formula  for almost all $x \in \R$, see Proposition \ref{p4}. The real analytic functions $\Theta^{\pm}$, $\Phi^{\pm}$ satisfy
$$
\Theta^+(r,z)\Phi^-(r,z) - \Theta^-(r,z)\Phi^+(r,z)  = \det M_0(r, z) = 1 
$$
for all $r \ge 0$, $z \in \C$, hence we obtain 
\begin{equation}\label{eq11}
w_0(x) = \Im m_0(x) = \frac{\Im m_r(x)}{|F_r(x)|^{2}} = \frac{w_r(x)}{|F_r(x)|^{2}},
\end{equation}
for almost all $x \in \R$, where $F_r: z\mapsto \Theta^+(r,z) + m_r(z)\Theta^-(r,z)$ is the analytic function in $\C^+$ and $F_r(x)$, $x \in \R$, are the non-tangential boundary values of $F_r$. Denote the first column of the matrix-function $M$ in \eqref{eq65}  by $\Theta=\left(\!\begin{smallmatrix}\Theta^+\\ \Theta^-\end{smallmatrix}\!\right)$. Assume for a moment that $(0,r)$ is not an indivisible interval of type $\pi/2$ for $\Hh$. Then formula \eqref{eq61} implies that $\Theta^{-}(r,z) \neq 0$ for every $z \notin \R$, and, moreover, $\Im\frac{\Theta^{+}(r, z)}{\Theta^-(r, z)}> 0$ for $z \in \C^+$. Thus, the function $\log|F_r|$ can be represented in the form
$$
\log|F_r(z)| = \log|\Theta^{-}(r,z)| + \log\left|m_r(z) + \frac{\Theta^{+}(r, z)}{\Theta^-(r, z)}\right|, \qquad z \in \C^+.
$$ 
Since the functions $m_r$, $\frac{\Theta^{+}(r, \cdot)}{\Theta^-(r, \cdot)}$ have positive imaginary parts in~$\C^+$ and $\Theta^{-} \in \N(\C^+)$, we have $|\log|F_r(x)||\,dx \in \Pi(\R)$, and, moreover,
$$
\frac{1}{\pi}\int_{\R}\frac{\log |F_r(x)|}{1+x^2}\,dx = \log|F_r(i)| - \xi_\Hh(r),
$$
by Proposition \ref{p2} and Proposition \ref{p4}. In particular, the measure $\mu_r$ belongs to the Szeg\H{o} class~$\sz$. Taking logarithms in \eqref{eq11} and integrating with $\frac{1}{1+x^2}$, we obtain assertion~$(b)$:
\begin{equation}\label{eq13}
\J_\Hh(r) =\J_{\Hh}(0) - 2\xi_\Hh(r) + 2\log |F_r(i)|. 
\end{equation}
Let us now prove $(b)$ in the case where $\Hh$ has an indivisible interval $(0, \eps)$ of type $\pi/2$ for some $\eps > 0$ and $r\le \eps$. In that situation, we can use Lemma \ref{l5} to show that $F_r(z) = 1$ for all $z$, hence $w_0 = w_r$ on $\R$ by \eqref{eq11}, yielding $\J_{\Hh}(r) = \J_{\Hh}(0)$ for $r\in [0,\eps]$. Since $\xi_{\Hh} = 0$ on $[0, \eps]$ by definition, this gives us relation $(b)$ in full generality.

\medskip
  
Next, the solution $M^d(r,z)$ of the canonical Hamiltonian system generated by the dual Hamiltonian $\Hh^d=J^\ast \Hh J$ has the form 
\begin{equation}\label{eq69}
M^d(r,z) = J^* M(r,z) J = \begin{pmatrix}\Phi^-(r,z) &  -\Theta^-(r,z) \\ -\Phi^+(r,z) & \Theta^+(r,z)\end{pmatrix}.
\end{equation}
Note that $\Hh^d$, $\Hh^d_r$ are singular nontrivial Hamiltonians because $\Hh$, $\Hh_r$ are singular and nontrivial. Using formula \eqref{eq9} with $\omega=\infty$, we see that $m_r^d(z) = -\lim_{t\to+\infty}\frac{\Theta_r^+(t,z)}{\Phi_r^+(t,z)} = -\frac{1}{m_r(z)}$ for all $r \ge 0$ and all $z \in \C^+$. Taking the non-tangential values of imaginary parts gives $w_{r}^{d}(x) = \frac{\Im m_r(x)}{|m_r(x)|^2} = \frac{w_r(x)}{|m_r(x)|^2}$. This formula  and Proposition \ref{p4} imply $\mu^d_r \in \sz$ thus completing the proof of $(a)$. Since the measures $\mu_r$, $\mu_r^d$ are even, we have 
\begin{equation}\label{eq23}
\I_{\Hh^d}(r) = \Im m_r^d(i) = \frac{1}{\Im m_r(i)} = \frac{1}{\I_{\Hh}(r)},
\end{equation} 
as claimed in $(c)$. Next, using the formula $w_{r}^{d}(x) = \frac{w_r(x)}{|m_r(x)|^2}, x\in \mathbb{R}$, the mean value formula in Proposition \ref{p4},  formula \eqref{eq23}, and identity $m_r(i) = i\I_{\Hh}(r)$, we obtain assertion $(d)$:
\begin{align*}
\K_{\Hh^d}(r) 
&= \log\I_{\Hh^d}(r) - \J_{\Hh}(r) + \log|m_r(i)|^2\\ 
&= -\log \I_{\Hh}(r) - \J_{\Hh}(r) + 2\log\I_{\Hh}(r) = \K_{\Hh}(r).
\end{align*}
Finally, consider the Hamiltonian $\widehat \Hh_r$ introduced in \eqref{cook59}.
 Since $\Hh_r$ is nontrivial, we have $\I_{\Hh}(r) \neq 0$ and hence $\widehat\Hh_r$ is defined correctly. By definition and Lemma \ref{l9}, we have $\I_{\widehat \Hh_r}(r) = \I_{\Hh}(r)$, $\J_{\widehat \Hh_r}(r) = \log\I_{\Hh}(r)$, and $\widehat F_r(i) = F_r(i)$ for the corresponding function $\widehat F_r$. The proof of Lemma \ref{l9} shows that $\widehat m_t$ is a constant function for each $t \ge r$. Using this and the fact that $\Phi^{\pm}, \Theta^{\pm} \in \N(\C^+)$, from \eqref{eq10} we obtain $\widehat\mu_r \in \sz$. Comparing the right hand sides of formula \eqref{eq10} for $m_0$ and $\widehat m_0$ at $z =i$, we get $\I_{\widehat \Hh_r}(0) = \I_{\Hh}(0)$. Hence, relation \eqref{eq13} for $\widehat H_r$ can be written in the form
$$
\J_{\widehat \Hh_r}(r) = \J_{\widehat\Hh_r}(0) - 2\xi_\Hh(r) + 2\log |F_r(i)| = \J_{\widehat\Hh_r}(0) - \J_{\Hh}(0) + \J_{\Hh}(r). 
$$
On the other hand, we have $\log\I_{\Hh}(r) = \J_{\widehat \Hh_r}(r)$ and $\I_{\widehat \Hh_r}(0) = \I_{\Hh}(0)$. This yields assertion $(e)$: 
\begin{align*}
\K_\Hh(r) 
&=\log\I_{\Hh}(r) - \J_{\Hh}(r) =\J_{\widehat \Hh_r}(r) - \J_{\Hh}(r) =\J_{\widehat\Hh_r}(0) -  \J_{\Hh}(0)\\
&= \J_{\widehat\Hh_r}(0) - \log\I_{\Hh}(0) + \log\I_{\Hh}(0) - \J_{\Hh}(0)\\
&= \J_{\widehat\Hh_r}(0) - \log\I_{\widehat \Hh_r}(0) + \log\I_{\Hh}(0) - \J_{\Hh}(0)\\
&= -\K_{\widehat\Hh_r}(0) + \K_{\Hh}(0). 
\end{align*} 
The lemma is proved. \qed
\medskip

\begin{Lem}\label{l0}
Let $l>0$ and $\Hh$ be a singular Hamiltonian on $\R_+$  satisfying   $\Hh(t) = \diag(a_1, a_2)$ for  all $t \in [\ell, + \infty)$ where $a_1$, $a_2$ are positive parameters. Then its spectral measure $\mu$ belongs to the Szeg\H{o} class $\sz$.
\end{Lem}
\beginpf Formula \eqref{eq11} for $r =\ell$ says that the absolutely continuous part of $\mu$ coincides with $\frac{|w_\ell(x)|^2}{|F_{\ell}(x)|^2}$. Since $\Hh_\ell = \diag(a_1, a_2)$ on $\R_+$,  we have $w_\ell(x) = \sqrt{a_2/a_1}$ for all $x \in \R$  by Lemma \ref{l9}. It remains to use Proposition \ref{p2} for the function $F_\ell \neq 0$ of class $\N(\C^+)$. \qed

\medskip

\begin{Lem}\label{l2}
Let $\Hh = \diag(h_1, h_2)$ be a singular nontrivial Hamiltonian on $\R_+$ whose spectral measure belongs to the Szeg\H{o} class $\sz$. Then the functions $\J_{\Hh}(r), \K_{\Hh}(r)$ are absolutely continuous and 
\begin{align}
\J'_{\Hh}(r) &= 2\I_{\Hh}(r)h_1(r) - 2\xi'_{\Hh}(r),\label{eq14} \\
\K'_{\Hh}(r) &= -\I_{\Hh}(r)h_1(r) -\frac{h_2(r)}{\I_{\Hh}(r)} + 2\xi'_{\Hh}(r),\label{eq54}
\end{align} 
for almost all $r \ge 0$.
\end{Lem} 
\beginpf At first, assume additionally that $h_1$, $h_2$ belong to $\Cc^1(\R_+)$, the space of continuously differentiable functions on $(0, +\infty)$ whose derivatives have a finite limit at~$0$. Then the entries of the the solution $M(\cdot,i)$ of \eqref{eq1} at $z=i$ belong to the space $\Cc^1(\R_+)$ as well. From formula \eqref{eq10} and identity $m_r(i) = i\I_{\Hh}(r)$, $r \ge 0$, we also have $\I_{\Hh} \in \Cc^1(\R_+)$. Assertion $(b)$ of Lemma~\ref{l1} says that
\begin{equation}\label{eq15}
\J_\Hh(r)= \J_{\Hh}(0) - 2\xi_\Hh(r) + 2\log|\Theta^+(r,i) +i\I_\Hh(r)\Theta^-(r,i)|, \qquad r \ge 0.
\end{equation}
Differentiating the above formula with respect to $r$ at $r = 0$ and using the equation 
$$
\left.\left(\!\begin{smallmatrix}\Theta^+(r,i)'&\Phi^+(r,i)'\\ \Theta^-(r,i)'&\Phi^-(r,i)'\end{smallmatrix}\!\right)\right|_{r=0} 
= M'(0,i) = i J^* \Hh(0) M(0,i) = \left(\!\begin{smallmatrix}0&ih_2(0)\\ -ih_1(0)&0\end{smallmatrix}\!\right),
$$ 
we obtain 
\begin{align*}
\J'_\Hh(0) 
&= -2\xi'_\Hh(0) + 2 \left.\Re\left(\frac{\Theta^+(r,i)' +i\I'_\Hh(r)\Theta^-(r,i) + i\I_\Hh(r)\Theta^-(r,i)'}{\Theta^+(r,i) +i\I_\Hh(r)\Theta^-(r,i)}\right)\right|_{r=0}\\
&= -2\xi'_\Hh(0) + 2\I_\Hh(0) h_1(0).
\end{align*}
For $r > 0$ we have
$$
\J'_\Hh(r) = \J'_{\Hh_r}(0) = -2\xi'_{\Hh_r}(0) + 2\I_{\Hh_r}(0) h_1(r) = -2\xi'_\Hh(r) + 2\I_{\Hh}(r) h_1(r).   
$$
Thus, relation \eqref{eq14} holds in the case when $h_1, h_2\in \Cc^1(\R_+)$.  Now let $\Hh = \diag(h_1,h_2)$ be an arbitrary singular nontrivial Hamiltonian on $\R_+$ with spectral measure in $\sz$. By Lemma \ref{l1}, the functions $\I_{\Hh}(r)$, $\J_\Hh(r)$  are correctly defined on $\R_+$. Find a sequence of positive smooth functions 
$\{h_{1,n}\}$, $\{h_{2,n}\}$ such that 
\[
\lim_{n\to \infty} \int_{0}^{T}|h_j(s) - h_{j,n}(s)|\,ds=0
 \]
for every $T>0$ and $j=1,2$. Solutions of the equations $JM'_{(n)} = i\Hh_{(n)} M_{(n)}$, $M_{(n)}(0,i) = \idm$, generated by the Hamiltonians $\Hh_{(n)} = \diag(h_{1,n}, h_{2,n})$ will then converge uniformly on compact subsets of $\R_+$ to the solution $M(\cdot, i)$ of the equation $JM' = i\Hh M$, $M(0,i) = \idm$. From formulas \eqref{eq10} and \eqref{eq15} we see that continuous functions $\I_{\Hh_{(n)}}(r)$, $\J_{\Hh_{(n)}}(r)$ converge  uniformly on compact subsets of $\R_+$ to the functions $\I_{\Hh}(r)$, $\J_{\Hh}(r)$, respectively. Thus, we have 
\begin{align*}
\J_\Hh(r) - \J_\Hh(0) 
&= \lim_{n \to \infty} (\J_{\Hh_{(n)}}(r) - \J_{\Hh_{(n)}}(0)) \\
&= -2 \xi_\Hh(r) + \lim_{n \to \infty}\int_{0}^{r}\I_{\Hh_{(n)}}(s) h_1(s)\,ds \\
&= -2 \xi_\Hh(r) + \int_{0}^{r}\I_{\Hh}(s) h_1(s)\,ds,
\end{align*}
for every $r > 0$. This formula shows that $\J_\Hh$ is absolutely continuous and satisfies relation \eqref{eq14}. Relation \eqref{eq54} follows by adding  \eqref{eq14} written for $\Hh$ and $\Hh_d = \diag(h_2, h_1)$  and using identity
\begin{equation}\label{sdcif}
\K_{\Hh}=-(\J_{\Hh}+\J_{\Hh_d})/2
\end{equation}
which is immediate from Lemma \ref{l1}.$(c),(d)$. \qed

\medskip
 
\begin{Lem}\label{l3}  Let $l>0$ and $\Hh$ be a singular Hamiltonian on $\R_+$  satisfying   $\Hh(t) = \diag(a_1, a_2)$ for  all $t \in [\ell, + \infty)$ where $a_1$, $a_2$ are positive parameters. 
Then, for every $r \ge 0$ we have
\begin{align}
e^{-\frac{1}{2}\J_{\Hh}(r) - \xi_{\Hh}(r)} &= \int_{r}^{\infty}h_1(s)e^{-\frac{1}{2}\J_{\Hh^d}(s)-\xi_\Hh(s)}\,ds, \label{eq16}\\ 
e^{-\frac{1}{2}\J_{\Hh^d}(r) - \xi_{\Hh}(r)} &= \int_{r}^{\infty}h_2(s)e^{-\frac{1}{2}\J_{\Hh}(s)-\xi_\Hh(s)}\,ds \label{eq17}.
\end{align}
\end{Lem}
\beginpf The right hand side of \eqref{eq16} at $r_0 \ge \ell$ is equal to 
$$
a_1 e^{-\xi_\Hh(r_0) - \tfrac{1}{2}\J_{\Hh^d}(r_0)}\int_{r_0}^{\infty}e^{(r_0-s)\sqrt{a_1 a_2}}\,ds 
= \sqrt{\frac{a_1}{a_2}} \cdot e^{-\xi_\Hh(r_0) - \tfrac{1}{2}\J_{\Hh^d}(r_0)}.
$$ 
Substituting 
$\J_{\Hh}(r_0) = \log\sqrt{\frac{a_2}{a_1}}$, $\J_{\Hh^d}(r_0) = \log\sqrt{\frac{a_1}{a_2}}$ into the formula above, we see that \eqref{eq16} holds for all $r \ge \ell$. Next, differentiating the left hand side of \eqref{eq16} and using Lemma~\ref{l1} and Lemma~\ref{l2}, we obtain
\begin{align*}
-\left(\!\frac{\J_{\Hh}'(r)}{2} + \xi'_\Hh(r)\!\right)\!e^{-\frac{1}{2}\J_{\Hh}(r) - \xi_{\Hh}(r)} 
&= -h_1(r) \I_{\Hh}(r) e^{-\frac{1}{2}\J_{\Hh}(r) - \xi_{\Hh}(r)}\\
&=-h_1(r) e^{\frac{1}{2}\log \I_{\Hh}(r) + \frac{1}{2} \K_{\Hh}(r) - \xi_{\Hh}(r)}\\
&=-h_1(r) e^{\frac{1}{2}\log \I_{\Hh}(r) + \frac{1}{2} (\log \I_{\Hh^d}(r) - \J_{\Hh^d}(r)) - \xi_{\Hh}(r)}\\
&=-h_1(r) e^{- \frac{1}{2}\J_{\Hh^d}(r) - \xi_{\Hh}(r)}.
\end{align*}
This agrees with the derivative of the right hand side of \eqref{eq16} for almost all $r \ge 0$.  It follows that \eqref{eq16} holds for all $r \ge 0$. Formula \eqref{eq17} can be proved in a similar way. \qed   

\medskip

\section{Some estimates of the entropy function}\label{s3}
In this section we consider Hamiltonians $\Hh$ such that $\det\Hh = 1$ almost everywhere on $\R_+$. In the notations of Section \ref{s2}, we have $\K(\mu) = \K_{\Hh}(0)$ for such Hamiltonians. Indeed, the coefficient $b_0$ in \eqref{eq63} is non-zero if and only if there exists $\eps>0$ such that $(0, \eps)$ is the indivisible interval of type $\pi/2$ for $\Hh_0 = \Hh$, see Lemma \ref{l12}. The latter never happens for Hamiltonians $\Hh$ with $\det\Hh = 1$ almost everywhere on $\R_+$.  
\subsection{A lower bound for the entropy}
We first obtain a local estimate for the entropy $\K(\mu) = \K_{\Hh}(0)$ in terms of $\Hh$ and then use assertion $(e)$ of Lemma \ref{l1} to improve it. 
\begin{Lem}\label{l4}
Let $h \ge 0$ be a function on $\R_+$ such that  $h, 1/h \in L^1_{\rm loc}(\R_+)$ and assume that $h$ equals to some positive constant on $[\ell, +\infty)$ for some $\ell\ge 0$. Then, for the Hamiltonian $\Hh = \diag(h,1/h)$, we have
$$
e^{\frac{1}{2}\K_{\Hh}(0)}\ge \int_{0}^{\infty}\sqrt{a(t)} \cdot te^{-t}dt, 
$$
where $a(t) = \frac{1}{t}\int_{0}^{t}h(s)\,ds \cdot \frac{1}{t}\int_{0}^{t}\frac{1}{h(s)}\,ds$ for $t>0$.
\end{Lem}
\beginpf Using Lemma \ref{l3} twice, we get
\begin{align}
e^{-\frac{1}{2}\J_{\Hh}(0)} 
&= \int_{0}^{\infty}h(s)e^{-\frac{1}{2}\J_{\Hh^d}(s) - s}\,ds \notag\\
&= \int_{0}^{\infty}h(s)\left(\int_{s}^{\infty}\frac{1}{h(\tau)}e^{-\frac{1}{2}\J_{\Hh}(\tau)}e^{s-\tau}\,d\tau\right) e^{-s}\,ds \notag\\
&= \int_{0}^{\infty}\frac{1}{h(\tau)}e^{-\frac{1}{2}\J_{\Hh}(\tau)}\left(\int_{0}^{\tau}h(s)\,ds\right) e^{-\tau}\,d\tau \label{eq18}.
\end{align}
Analogous formula holds for $\J_{\Hh^d}$:
\begin{equation}\label{eq19}
e^{-\frac{1}{2}\J_{\Hh^d}(0)} 
= \int_{0}^{\infty}h(\tau)e^{-\frac{1}{2}\J_{\Hh^d}(\tau)}\left(\int_{0}^{\tau}\frac{1}{h(s)}\,ds\right) e^{-\tau}\,d\tau.
\end{equation}
We have $2\K_{\Hh}(r) = -\J_{\Hh}(r) - \J_{\Hh^d}(r)$ for all $r \ge 0$ (see \eqref{sdcif}). We also have $\K_\Hh \ge 0$ on $\R_+$  (check, e.g.,   \eqref{sdp1}). Multiplying formulas \eqref{eq18}, \eqref{eq19} and using Cauchy-Schwarz inequality, we obtain
$$
e^{\frac{1}{2}\K_{\Hh}(0)} 
\ge \int_{0}^{\infty}e^{\frac{1}{2}\K_\Hh(\tau)}e^{-\tau}\sqrt{\int_{0}^{\tau}h(s)\,ds\int_{0}^{\tau}\frac{1}{h(s)}\,ds}\,d\tau
\ge \int_{0}^{\infty}\sqrt{a(t)} \cdot te^{-t}dt,
$$
as required. \qed\medskip

\medskip

\noindent{\bf Remark.} We can write $a(t)=\langle h\rangle_{[0,t]}\langle 1/h\rangle_{[0,t]}$ and  $a(t)\ge 1$, as follows from  Cauchy-Schwarz inequality.

\medskip

This lemma and additivity of the entropy $\K_{\Hh}$ imply the following estimate. 
\begin{Prop}\label{p1}
Let $h \ge 0$ be a function on $\R_+$ such that  $h, 1/h \in L^1_{\rm loc}(\R_+)$ and $\Hh = \diag(h, 1/h)$. Then, there exists a sequence of numbers $\{t_n\}$ such that $t_n \in [3, 4]$ and
$$
\sum_{n \ge 0} \left(\frac{1}{t_{n}}\int_{4n}^{4n+t_{n}}h(s)\,ds \cdot \frac{1}{t_{n}}\int_{4n}^{4n+t_{n}}\frac{ds}{h(s)} -1 \right) \le e^{10\K_{\Hh}(0)} - 1.
$$ 
\end{Prop}
\beginpf Iteratively applying assertion $(e)$ of Lemma \ref{l1}, we can find a sequence of Hamiltonians $\Hh_{(n)} = \diag(h_{n},1/h_{n})$ such that
$\Hh_{(n)}(x) = \Hh(4n+x)$ for $x \in [0, 4]$, $\Hh_{(n)}(x) = \diag(a_{n},1/a_{n})$ for almost all $x > 4$ and some constant $a_n>0$, and 
\begin{equation}\label{sdp2}
\K_{\Hh}(0) \ge  \sum_{n \ge 0}\K_{\Hh_{(n)}}(0).  
\end{equation}
Take $n \ge 0$ and apply Lemma \ref{l4} for the Hamiltonian~$\Hh_{(n)}$. Making note of 
\[
\int_0^{\infty} te^{-t}dt=1
\]
and applying  Jensen inequality, we get
$$
\K_{\Hh_{(n)}}(0) \ge \int_{0}^{\infty}\log a_n(t) \cdot t e^{-t}\,dt,
$$
where $a_n(t) = \frac{1}{t}\int_{4n}^{4n+t}h(s)\,ds \cdot \frac{1}{t}\int_{4n}^{4n+t}\frac{1}{h(s)}\,ds$ for $t \in [0,4]$ and $a_n(t) \ge 1$ for all $t>0$. Since $\int_{I}te^{-t}dt \ge 0.1$ for $I = [3, 4]$, we have $10\K_{\Hh_{(n)}}(0) \ge \min_{t \in I}\log a_n(t)$. Define $t_n$ to be a point in $I$ such that $a_n(t_n) = \min_{t \in I}a_n(t)$. Since $e^{x+y} - 1 \ge e^{x} - 1 + e^y - 1$ for all $x, y \ge 0$, we notice that \eqref{sdp2} implies
\begin{align*}
e^{10\K_{\Hh}(0)} - 1
&\ge \sum_{n \ge 0} \left(e^{10\K_{\Hh_{(n)}}(0)} - 1\right) \ge \sum_{n \ge 0} \left(a_n(t_n) - 1\right)\\
&= \sum_{n \ge 0} \left(\frac{1}{t_n}\int_{4n}^{4n+t_n}h(s)\,ds \cdot \frac{1}{t_n}\int_{4n}^{4n+t_n}\frac{1}{h(s)}\,ds - 1\right),
\end{align*}   
which is the desired estimate. \qed

\medskip

\subsection{An upper bound for the entropy}
\begin{Prop}\label{p141}
Let $h$ be a function as in Lemma \ref{l4}, and let $\Hh = \diag(h,1/h)$ be the corresponding Hamiltonian. Then, 
$$
\K_{\Hh}(0) \le \int_{0}^{\infty}(\kappa(s) + \kappa_d(s) - 2)\,ds,
$$
where $\kappa(r) = \frac{1}{h(r)}\int_{r}^{\infty}h(s)e^{r-s}\,ds$ and $\kappa_d(r) = h(r)\int_{r}^{\infty}\frac{1}{h(s)}e^{r-s}\,ds$ for $r \ge 0$.
\end{Prop}
\beginpf Consider the functions  
$$
u(r) = \int_{r}^{\infty}\frac{1}{h(s)}e^{-\J_{\Hh}(s)-s}\,ds,  \qquad u_d(r) = \int_{r}^{\infty}h(s)e^{-\J_{\Hh^d}(s)-s}\,ds,
$$
defined on $\R_+$. By Lemma \ref{l3}, we have 
\begin{align*}
e^{-\J_{\Hh}(r)} 
&= \left(\int_{r}^{\infty}h(s)e^{-\frac{\J_{\Hh^d}(s)}{2}}e^{r-s}\,ds\right)^2\\ 
&\le \left(\int_{r}^{\infty}h(s)e^{r-s}\,ds\right) \left(\int_{r}^{\infty}h(s)e^{-\J_{\Hh^d}(s)}e^{r-s}\,ds\right)\\
&= h(r)e^r \kappa(r)u_d(r).
\end{align*}
Dividing by $he^{r}$, we obtain $-u'(r) \le \kappa(r) u_d(r)$ for almost all $r \ge 0$. Analogously, we have $-u'_d(r) \le \kappa_d(r) u(r)$, $r \ge 0$ for the function $u_d$. It follows that
$$
0 \le -(u^2 + u_d^2)'(r) \le 2(\kappa(r) + \kappa_d(r))u(r) u_d(r) \le (\kappa(r) + \kappa_d(r))(u^2 + u_d^2)(r),
$$
for almost all $r \ge 0$. Thus, we have
$$
-\frac{\partial}{\partial r}\log\bigl(u^2(r) + u_d^2(r)\bigr) \le \kappa(r) + \kappa_d(r).
$$
Taking into account that $u(r) = u_d(r) = e^{-r}$ for $r \ge \ell$ by \eqref{sdp3}, we get
\begin{equation}\label{eq20}
u^2(0) + u_d^2(0) \le (u^2(\ell) + u_d^2(\ell)) e^{\int_{0}^{\ell}(\kappa(s) + \kappa_d(s))\,ds} = 2e^{\int_{0}^{+\infty}(\kappa(s) + \kappa_d(s) - 2)\,ds}\,.
\end{equation}
On the other hand, we have 
$$u(0) = \int_{0}^{\infty}\frac{1}{\I_{\Hh}(s)h(s)}e^{\K_{\Hh}(s) - s}\,ds, \qquad 
u_d(0) = \int_{0}^{\infty}\I_{\Hh}(s)h(s)e^{\K_{\Hh}(s) - s}\,ds,$$ 
by assertions $(c)$, $(d)$ of Lemma \ref{l1}. From \eqref{eq54} for $h_1 = h = 1/h_2$ we now get
\begin{align*}
u(0) + u_d(0) 
&= - \int_{0}^{\infty}\K'_{\Hh}(s)e^{\K_{\Hh}(s)-s}\,ds + 2\int_{0}^{\infty}e^{\K_{\Hh}(s)-s}\,ds\\ 
&= e^{\K_{\Hh}(0)} + \int_{0}^{\infty}e^{\K_{\Hh}(s)-s}\,ds\\
&\ge e^{\K_{\Hh}(0)} + 1 \ge 2e^{\K_{\Hh}(0)/2}, 
\end{align*}
using integration by parts and  the fact that $\K_{\Hh}(s) \ge 0$ for all $s$. Last estimate and \eqref{eq20} imply
$$
e^{\K_{\Hh}(0)} \le \left(\frac{u(0) + u_d(0)}{2}\right)^2 \le \frac{u^2(0) + u^2_d(0)}{2} \le e^{\int_{0}^{+\infty}(\kappa(s) + \kappa_d(s) - 2)\,ds}.  
$$
Taking the logarithms, we arrive to the statement of the proposition. \qed

\medskip


\section{Proof of Theorem \ref{t1}}\label{s4}
The classical Muckenhoupt  class $A_2(\R)$ is defined as the set of measurable functions $h \ge 0$ on $\R$ with finite characteristic 
\[
[h]_{2}\equiv \sup_{I \subset \R}\,\langle h\rangle_I\langle h^{-1}\rangle_I,
\]
where the supremum is taken over all intervals $I \subset \mathbb{R}$. Recall that $I_{x,y}$ denotes $[x,x+y)$ for $x, y \in \R_+$. For a function $h \ge 0$ on $\R_+$ and a sequence $\alpha = \{\alpha_n\}$ of positive numbers, put  
\begin{equation}\label{eq30}
[h, \alpha] = \sum_{n=0}^\infty \Bigl(\langle h\rangle_{I_{n,\alpha_n}} \langle h^{-1}\rangle_{{I_{n,\alpha_n}}}-1\Bigr).
\end{equation} 
Each term in the sum above is nonnegative, hence $[h, \alpha] \in \R_+ \cup \{+\infty\}$ is correctly defined. Denote by $\two$ the constant sequence $2,2, \ldots$ indexed by non-negative integers. 

\medskip

\noindent{\bf Definition.} Let $A_2(\R_+, \ell^1)$ be the set of functions $h \ge 0$ on $\R_+$ such that the characteristic $[h]_{2,\, \ell^1} = [h, \two]$ is finite.

\medskip

Note that $[h]_{2,\, \ell^1} = 0$ if and only if the function $h$ is constant. Next, for a function $h \ge 0$ on $\R_+$ define 
\begin{equation}\label{eq29}
[h]_{int} = \int_{0}^{\infty}(\kappa(s) + \kappa_d(s) - 2)\,ds, 
\end{equation}
where $\kappa(r) = \frac{1}{h(r)}\int_{r}^{\infty}h(s)e^{r-s}\,ds$ and $\kappa_d(r) = h(r)\int_{r}^{\infty}\frac{1}{h(s)}e^{r-s}\,ds$ for $r \ge 0$. Since $h \ge 0$ on $\R_+$,  we have $\frac{h(s)}{h(r)} + \frac{h(r)}{h(s)} \ge 2$, hence the quantity $[h]_{int} \in \R_+ \cup \{+\infty\}$ is correctly defined.  
\begin{Prop}\label{p9}
Let $h \ge 0$ be a measurable function on $\R_+$. Assume that 
$[h, \alpha]$ is finite for a sequence $\alpha=\{\alpha_n\}$ where $\alpha_n\in  [3, 4], \forall n\in \mathbb{Z}^+$. Then $h \in A_2(\R_+, \ell^1)$ and, moreover, we have  
$[h]_{2,\, \ell^1} \le c[h, \alpha]$ with absolute constant $c$.
\end{Prop}
\begin{Prop}\label{p5}
There exists an absolute constant $c$ such that  $[h]_{int} \le c[h]_{2,\ell^1}e^{c[h]_{2,\ell^1}}$ for every function $h \in A_2(\R_+, \ell^1)$ . 
\end{Prop}
Propositions \ref{p9}, \ref{p5} will be proved in the next section. Later, in the proof of the theorem, we will need the following lemma.
\begin{Lem}\label{l6}
Let $\Hh$, $\Hh_{(k)}$ be singular diagonal Hamiltonians on $\R_+$ such that $\Hh_{(k)}(x) = \Hh(x)$ for every $k \ge 0$ and all $x \in [0,k]$. 
Suppose that the spectral measure of $\Hh_{(k)}$ belongs to $\sz$  for every $k\ge 0$  and $\sup_{k \ge 0}\K_{\Hh_{(k)}}(0) < \infty$. Then, the spectral measure of $\Hh$ belongs to $\sz$  and $\K_{\Hh}(0) \le \limsup_{k \to \infty}\K_{\Hh_{(k)}}(0)$. 
\end{Lem}
\beginpf Let $\Hh$ be a singular Hamiltonian on $\R_+$ and let $m$ be its Weyl-Titchmarsh function. As usual, denote by $\Theta^{\pm}$, $\Phi^{\pm}$ the corresponding entries of the solution $M$ of Cauchy problem \eqref{eq1}. Then, by the nesting circles analysis (see page 42 in Section 8 of \cite{Romanov} or page 475 in Section 7 of \cite{HSW}), we have
\begin{equation}\label{eq24}
\left|m(z) - \frac{\Phi^{-}(k,z)}{\Theta^{-}(k,z)}\right| \le \frac{1}{\Im\bigl(\Theta^+(k,z)\ov{\Theta^-(k,z)}\bigr)} , \qquad z \in \C^+, \quad k \ge 0,
\end{equation}
where the right hand side tends to zero as $k \to +\infty$ uniformly on compacts in $\C^+$. Let $m_{(k)}$ be the Weyl-Titchmarsh function of the Hamiltonian $\Hh_{(k)}$. Since $\Hh_{(k)}$ coincides with $\Hh$ on $[0,k]$, we have estimate \eqref{eq24} with $m$ replaced by $m_{(k)}$ and the same right hand side. The triangle inequality now implies that $m - m_{(k)}$ tends to zero uniformly on compact subsets of $\C^+$. 

\medskip

Let us consider the measures $\widetilde\mu$, $\widetilde\mu_{(k)}$ supported on the unit circle $\T = \{z \in \C: \; |z| = 1\}$ whose Poisson extensions to the open unit disk $\D = \{z \in \C: \; |z| < 1\}$ coincide with positive harmonic functions $\Im m(\omega)$, $\Im m_{(r)}(\omega)$ in $\D$, respectively, where $\omega: w \mapsto i\frac{1-w}{1+w}$ is the conformal mapping from $\D$ onto $\C^+$. Since the difference $m - m_{(k)}$ tends to zero uniformly on compacts in $\C^+$, the measures $\widetilde\mu_{(k)}$ converge weakly to the measure $\widetilde\mu$. Recall that the the relative entropy of two positive finite measures $\nu_1$, $\nu_2$ on $\T$ is defined by
$$
S(\nu_1|\nu_2) = \begin{cases}-\infty &\mbox{if $\nu_1$ is not $\nu_2$ a.c.},\\ - \int_{\T}\log\left(\frac{d\nu_1}{d\nu_2}\right)\,d\nu_1&\mbox{if $\nu_1$ is $\nu_2$ a.c.}.\end{cases}
$$ 
It is known (see Section 2.2.3 in \cite{Simonbook}) that the relative entropy is weakly upper-semicontinuous, which means $\limsup_{k \to +\infty} S(\nu_{1}|\nu_{2,k}) \le S(\nu_{1}|\nu_2)$ for every sequence of finite measures $\nu_{2, k}$ on $\T$  converging weakly to a measure $\nu_2$. 
This implies that $\widetilde\mu$ belongs to the Szeg\H{o} class on $\T$ and 
\begin{equation}\label{eq31}
- \infty < \limsup_{k \to \infty}\int_{\T}\log \widetilde w_{(k)}(\xi)\,dm(\xi) \le \int_{\T}\log \widetilde w(\xi)\,dm(\xi),
\end{equation} 
where $m$ is the Lebesgue measure on $\T$ normalized by $m(\T) = 1$, and $\widetilde w$, $\widetilde w_{(k)}$ are the densities on $\widetilde\mu$, $\widetilde\mu_{(k)}$ with respect to $m$. Changing variables  in \eqref{eq31}, we see that the spectral measure of $\Hh$ lies in the class $\sz$, and, moreover, 
$$
\limsup_{k\to+\infty} \J_{\Hh_{(k)}}(0) \le \J_{\Hh}(0). 
$$
From the relation $\lim_{k \to \infty} m_{(k)}(i) = m(i)$ we get $\I_{\Hh}(0) = \lim_{k \to +\infty}\I_{\Hh_{(k)}}(0)$. The lemma now follows. \qed 

\medskip
 
The next result establishes the key two-sided estimates for a special class of Hamiltonians.  
\begin{Lem}\label{l8}
Let $h$ be a function as in Lemma \ref{l4}, and let $\Hh = \diag(h,1/h)$. Then, we have $\K_\Hh(0) \le c\widetilde\K(\Hh) e^{c\widetilde\K(\Hh)}$ and $\widetilde\K(\Hh) \le c\K_\Hh(0) e^{c\K_\Hh(0)}$ for an absolute constant $c$. 
\end{Lem}
\beginpf By Lemma~\ref{l0}, the spectral measure of $\Hh$ belongs to $\sz$. From Proposition \ref{p141} we know that $\K_\Hh(0)\le [h]_{int}$. Proposition \ref{p5} implies $[h]_{int} \le c  [h]_{2,\ell^1} e^{c  [h]_{2,\ell^1}}$ with $ [h]_{2,\ell^1} = \widetilde \K(\Hh)$. Combining these estimates, we obtain inequality $\K_\Hh(0) \le c\widetilde\K(\Hh) e^{c\widetilde\K(\Hh)}$. To prove the second inequality, observe that  Proposition~\ref{p1}, when applied to $\Hh$, provides a sequence $\{t_n\} \subset [3, 4]$ such that
$$
\sum_{n \ge 0} \left(\frac{1}{t_{n}}\int_{4n}^{4n+t_{n}}h(s)\,ds \cdot \frac{1}{t_{n}}\int_{4n}^{4n+t_{n}}\frac{ds}{h(s)} -1 \right) \le e^{10\K_{\Hh}(0)} - 1.
$$
The same proposition applied to three ``translated'' Hamiltonians $\Hh_k: x \mapsto \Hh(x+k), k=1,2,3$, gives
$$
\sum_{n \ge 0} \left(\frac{1}{ t^{(k)}_{n}}\int_{4n}^{4n+ t^{(k)}_{n}}h(s+k)\,ds \cdot \frac{1}{ t^{(k)}_{n}}\int_{4n}^{4n+ t^{(k)}_{n}}\frac{ds}{h(s+k)} -1 \right) \le e^{10\K_{\Hh_k}(0)} - 1.
$$
for three new sequences $\{ t^{(k)}_n\} \subset [3, 4]$ where $k=1,2,3$.
Summing up the above four formulas, we obtain $[h, \alpha] \le e^{10 \K_{\Hh}(0)}-1 + \sum_{k=1}^3(e^{10 \K_{\Hh_{k}}(0)}-1)$ for the sequence $\alpha = \{\alpha_n\}$ defined by $\alpha_{4n} = t_n$, $\alpha_{4n+k} = t^{(k)}_{n}$, $n \ge 0, k=1,2,3$. By Lemma \ref{l1}.$(e)$, we have $\K_{\Hh_k}(0) \le \K_{\Hh}(0)$, hence 
$[h, \alpha] \le 4(e^{10\K_{\Hh}}(0) -1) \le c \K_{\Hh}(0)e^{10\K_{\Hh}(0)}$. 
Proposition \ref{p9} says that $[h]_{2, \ell^1} \le c [h, \alpha]$ for an absolute constant $c$. By definition, we have $\widetilde\K(\Hh) = [h]_{2, \ell^1}$, hence $\widetilde\K(\Hh) \le c\K_{\Hh}(0) e^{10\K_{\Hh}(0)}$. \qed

\medskip
In the next lemma, we will show that the condition that the determinant equals to one can be dropped.

\begin{Lem}\label{l7}
Let $\Hh = \diag(h_1, h_2)$ be a singular Hamiltonian on $\R_+$ such that $h_1$, $h_2$ are equal to positive constants on $[\ell, +\infty)$ for some $\ell \ge 0$. Then, we have $\widetilde\K(\Hh) \le c\K_\Hh(0) e^{c\K_\Hh(0)}$ and $\K_\Hh(0) \le c\widetilde\K(\Hh) e^{c\widetilde\K(\Hh)}$ with an absolute constant $c$.
\end{Lem}
\beginpf For every $\eps>0$ define $\Hh_{(\eps)}: t \mapsto \Hh(t) + \eps \chi_{[0,\ell]}(t) I_2$, $t \in \R_+$, where $I_2=\idm$ is the $2\times2$ identity matrix and $\chi_{[0,\ell]}$ denotes the characteristic function of $[0,\ell]$. 
Set $\xi_\eps = \xi_{\Hh_{(\eps)}}$, and let $\eta_\eps$ denote the inverse function to $\xi_\eps$, so that $\eta_{\eps}(\xi_{\eps}(t)) = t$ for all $t \ge 0$. 
Since $\xi_{\Hh_{(\eps)}}$ maps $\R_+$ onto $\R_+$, the function $\eta_\eps$ is defined correctly. Moreover, we have $\det \Hh_{(\eps)} > 0$ almost everywhere on $\R_+$, hence $\eta_\eps$ is absolutely continuous on $\R_+$ and we can define the Hamiltonian $\widetilde \Hh_{(\eps)}: t \mapsto \eta'_{\eps}(t)\Hh_{(\eps)}(\eta_\eps(t))$. 
By construction, $\eta'_\eps(t) = 1/\sqrt{\det\Hh_{(\eps)}(\eta_\eps(t))}$ almost everywhere on $\R_+$, so the Hamiltonian $\widetilde \Hh_{(\eps)}$ has determinant equal to one almost everywhere on $\R_+$. By Lemma \ref{l0}, the spectral measures $\mu$, $\mu_{(\eps)}$, $\widetilde \mu_{(\eps)}$ of $\Hh$, $\Hh_{(\eps)}$, $\widetilde \Hh_{(\eps)}$, respectively, belong to $\sz$. By Lemma~\ref{l8}, 
\begin{equation}\label{eq25}
\widetilde\K(\widetilde \Hh_{(\eps)}) \le c\K_{\widetilde\Hh_{(\eps)}}(0) e^{c\K_{\widetilde\Hh_{(\eps)}}(0)}, \qquad 
\K_{\widetilde\Hh_{(\eps)}}(0) \le c\widetilde\K(\widetilde \Hh_{(\eps)}) e^{c\widetilde\K(\widetilde \Hh_{(\eps)})},
\end{equation} 
for an absolute constant $c$. Let $h_{1,\eps}$, $h_{2,\eps}$, $h_\eps$ be defined by $\Hh_{(\eps)} = \diag(h_{1,\eps}, h_{2,\eps})$, $\widetilde \Hh_{(\eps)} = \diag(h_{\eps}, 1/h_{\eps})$. Then, for every $t \ge 0$, we have
$$
\int_{\eta_{\eps}(t)}^{\eta_{\eps}(t+2)}h_{1,\eps}(s)\,ds \cdot \int_{\eta_{\eps}(t)}^{\eta_{\eps}(t+2)} h_{2,\eps}(s)\,ds 
= \int_{t}^{t+2}h_{\eps}(s)\,ds \cdot \int_{t}^{t+2}\frac{1}{h_\eps(s)}\,ds,
$$
by a change of variables. This shows that $\widetilde \K(\widetilde \Hh_{(\eps)}) = \widetilde \K(\Hh_{(\eps)})$. It is also not difficult to see that the spectral measures $\mu_{(\eps)}$, $\widetilde\mu_{(\eps)}$ of $\Hh_{(\eps)}$, $\widetilde \Hh_{(\eps)}$ coincide. Indeed, solutions $M_{(\eps)}$, $\widetilde M_{(\eps)}$ of Cauchy problem 
\eqref{eq1} for $\Hh_{(\eps)}$, $\widetilde \Hh_{(\eps)}$ satisfy $\widetilde M_{(\eps)}(x) = M_{(\eps)}(\eta_{\eps}(x))$, $x \in \R_+$. Hence the limit in the right hand side of \eqref{eq372} defines the same harmonic function for $\Hh_{(\eps)}$ and $\widetilde \Hh_{(\eps)}$. Thus, from \eqref{eq25} we get  
\begin{equation}\label{eq38}
\widetilde\K(\Hh_{(\eps)}) \le c\K_{\Hh_{(\eps)}}(0) e^{c\K_{\Hh_{(\eps)}}(0)}, \qquad 
\K_{\Hh_{(\eps)}}(0) \le c\widetilde\K(\Hh_{(\eps)}) e^{c\widetilde\K(\Hh_{(\eps)})},
\end{equation} 
for every $\eps>0$. Next, by construction, we have $\xi_{\Hh_{(\eps)}}(t) > \xi_{\Hh}(t)$ for all $t>0$ and  $\eps>0$. Moreover, the difference $\xi_{\Hh_{(\eps)}} - \xi_{\Hh}$ tends to zero uniformly on $\R_+$ as $\eps$ tends to zero. Hence $\eta_{\eps}(t) < \eta(t)$ for all $t>0$, $\eps>0$ and $\eta(t) - \eta_{\eps}(t)$ tends to zero for each $t \in \R_+$ as $\eps$ tends to zero. Since $\Hh$, $\Hh_{(\eps)}$ are constant on $[\ell, + \infty)$, we have
\begin{align*}
0 &= \int_{\eta_n}^{\eta_{n+2}}h_{1}(s)\,ds \cdot \int_{\eta_n}^{\eta_{n+2}}h_{2}(s)\,ds  - 4,\\
0 &= \int_{\eta_{\eps}(n)}^{\eta_{\eps}(n+2)}h_{1,\eps}(s)\,ds \cdot \int_{\eta_{\eps}(n)}^{\eta_{\eps}(n+2)} h_{2,\eps}(s)\,ds  - 4,
\end{align*}
for all $n \ge n_0$ and all sufficiently small $\eps>0$, where $n_0$ can be chosen independently of~$\eps$. Hence, the sums in \eqref{eq2} which define $\widetilde \K(\Hh)$, $\widetilde\K(\Hh_{(\eps)})$ contain at most $n_0$ nonzero terms for small $\eps>0$. It follows that $\lim_{\eps \to 0} \widetilde\K(\Hh_{(\eps)}) = \widetilde\K(\Hh)$. 
It remains to show that $\lim_{\eps \to 0}\K_{\Hh_{(\eps)}}(0) = \K_{\Hh}(0)$. To do that, one can use formula \eqref{eq10} with $r = \ell$ for $\Hh$ and $\Hh_{(\eps)}$. Since the matrix norm of $\Hh - \Hh_{(\eps)}$ tends to zero uniformly on $[0, \ell]$ and $\Hh = \Hh_{(\eps)}$ on $[\ell, + \infty)$, we have   
\begin{equation}\label{eq66}
\J_{\Hh}(\ell) = \J_{\Hh_{(\eps)}}(\ell), \quad \lim_{\eps \to 0}\xi_{\Hh_{(\eps)}}(\ell) = \xi_{\Hh}(\ell), \quad    
\lim_{\eps \to 0}|F_{\ell, \eps}(i)| = |F_{\ell}(i)|.
\end{equation}
To show that the last equality holds, we notice that  the Hamiltonians $\Hh_\ell$ and $\Hh_{(\eps)}(\cdot+\ell)$ coincide on $\R_+$ and thus have the same Weyl-Titchmarsh functions which we denote by $m_\ell$. Hence, the corresponding functions $F_{\ell, \eps}: z \mapsto \Theta^{+}_{(\eps)}(l,z) + m_\ell(z) \Theta^{-}_{(\eps)}(l,z)$ tend to $F_\ell$ uniformly on compact subsets of $\C^+$ as $\eps\to 0$. From \eqref{eq66} and Lemma \ref{l1}.$(b)$ for $r = \ell$,  we get $\lim_{\eps \to 0}\J_{\Hh_{(\eps)}}(0) = \J_{\Hh}(0)$. Using again formula \eqref{eq10} with $r = \ell$, we obtain $\lim_{\eps \to 0}\I_{\Hh_{(\eps)}}(0) = \I_{\Hh}(0)$. This completes the proof of the lemma. \qed 

\medskip
 Now we are ready to prove Theorem \ref{t1}.

\medskip

\noindent{\bf Proof of Theorem \ref{t1}.} Let $\Hh$ be a nontrivial singular diagonal Hamiltonian on $\R_+$ such that its spectral measure $\mu$ lies in the class $\sz$ and $b = 0$ in the Herglotz representation \eqref{eq371} of its Weyl-Tichmarsh function $m$. Note that we have $\K(\mu) = \K_{\Hh}(0)$ and  no positive $\eps$ exists such that $(0, \eps)$ is the indivisible interval for $\Hh$ of type $\pi/2$, see Lemma \ref{l12}. Consider the family of Bernstein-Szeg\H{o} Hamiltonians $\widehat \Hh_r = \diag(\widehat h_{1r}, \widehat h_{2r})$, $r \ge 0$, generated by $\Hh$ (see \eqref{cook59} for their definition).  By Lemma~\ref{l0},  the spectral measure $\widehat \mu_r$ of $\widehat \Hh_r$ belongs to $\sz$ for every $r \ge 0$. Since the Hamiltonians $\widehat \Hh_r$ have no indivisible intervals $(0, \eps)$ of type $\pi/2$, we have $\K(\widehat\mu_r) = \K_{\widehat\Hh_r}(0)$.  From Lemma~\ref{l1}.$(e)$ we now get $\K(\widehat \mu_r) \le \K(\mu)$. Let us first show that $\sqrt{\det\Hh} \notin L^1(\R_+)$. Since $2\sqrt{\det\Hh} \le \trace \Hh$, the function $\sqrt{\det\Hh}$ is integrable on compact subsets of $\R_+$. Suppose that $\sqrt{\det\Hh} \in L^{1}(\R_+)$. Then the function $\xi_\Hh$ in \eqref{eq7} is bounded, hence there exists $n_0 \ge 0$ and $r_0 \ge \eta_{n_0} \ge 0$, such that for every $r \ge r_0$ the last {\it nonzero} term in the sum defining $\widetilde\K(\widehat \Hh_r)$ equals 
$$
c_{r,n_0} = \int_{\eta_{n_0}}^{\widehat\eta_{n_0+2}(r)}\widehat h_{1r}(s)\,ds \cdot \int_{\eta_{n_0}}^{\widehat\eta_{n_0+2}(r)}\widehat h_{2r}(s)\,ds - 4,
$$  
where $\eta_{n_0} = \min\{t\ge 0: \xi_{\Hh}(t) = n_0\}$, and $\widehat\eta_{n_0+2}(r) = \min\{t\ge 0: \xi_{\widehat\Hh_r}(t) = n_0+2\}$ increases infinitely with $r$. By Lemma \ref{l7} and Lemma~\ref{l1}.$(e)$, we have $c_{r,n_0} \le \widetilde \K(\widehat\Hh_{r}) \le c\K(\widehat \mu_r)e^{c\K(\widehat \mu_r)} 
\le c\K(\mu)e^{c\K(\mu)}$ for every $r$.  From $\trace \Hh \notin L^1(\R_+)$ (recall that the Hamiltonian $\Hh$ is  singular) and the uniform boundedness of $c_{r, n_0}$, $r \ge r_0$, we get
\[
\int_{\eta_{n_0}}^\infty h_1(s)ds\int_{\eta_{n_0}}^\infty h_2(s)ds\le \limsup_{r\to\infty} c_{r,n_0}+4<\infty, \qquad
\int_{0}^\infty (h_1(s)+h_2(s))ds=\infty,
\]
which implies that either $\int_{\eta_{n_0}}^\infty h_1(s)ds=0$ or $\int_{\eta_{n_0}}^\infty h_2(s)ds=0$.  We see that ether $h_1 = 0$ or $h_2 = 0$ almost everywhere on $[r_0, +\infty)$ and  the Hamiltonian $\Hh_{r_0}$ is trivial. The first part of the proof of Lemma \ref{l1} shows that this is not the case, hence $\int_{0}^{\infty}\sqrt{\det\Hh(s)}\,ds = + \infty$\footnote{There is a different way to prove this fact. One needs to check that the supremum of the function $\xi_\Hh$ in \eqref{eq7} determines the exponential type of the measure $\mu$ and then apply Krein-Wiener completeness theorem. See Section 6 in \cite{Romanov}.} and the function $\eta_x$ in the statement of Theorem \ref{t1} is correctly defined on $\R_+$. For every $r \ge \eta_2$ the first $[\xi_\Hh(r)]-2$ terms defining $\widetilde \K(\Hh)$ and $\widetilde \K(\widehat \Hh_r)$ in \eqref{eq2} are identical. Hence, 
$$
\widetilde\K(\Hh) \le \limsup_{r \to \infty}\widetilde\K(\widehat \Hh_r) \le \limsup_{r \to \infty} c\K(\widehat \mu_r)e^{c\K(\widehat \mu_r)} 
\le c\K(\mu)e^{c\K(\mu)},
$$
where the second and the third inequalities follow from Lemma \ref{l7} and Lemma~\ref{l1}.$(e)$, respectively. 

\medskip

Conversely, suppose that $\Hh  = \diag(h_1, h_2)$ is a singular Hamiltonian on $\R_+$, $\sqrt{\det\Hh} \notin L^1(\R_+)$, and the sum defining $\widetilde\K(\Hh)$ in \eqref{eq2} converges. For every integer $k \ge 0$, fix some positive constants $a_{1k}$, $a_{2k}$ to be specified later, and consider
$$
\widetilde\Hh_{(k)}(t) = \diag(h_{1k}, h_{2k}) =
\begin{cases}
\Hh(t) & \mbox{ if } t \in [0,\eta_{k+2}],\\
\diag(a_{1k}, a_{2k})& \mbox{ if } t \in (\eta_{k+2}, +\infty).
\end{cases}
$$
For every $t>0$, set $\widetilde\eta_{t} = \min\{s \ge 0: \xi_{\Hh_{(k)}}(s) = t\}$, where $\xi_{\Hh_{(k)}}(s) = \int_{0}^{s}\sqrt{\det\Hh_{(k)}(\tau)}\,d\tau$. Then we have $\widetilde\eta_{t} = \eta_{t}$ for every $t \in [0,\eta_{k+2}]$. By construction, 
\begin{align}
\widetilde \K(\widetilde \Hh_{(k)}) 
=&\sum_{n = 0}^{k}\left(\int_{\eta_n}^{\eta_{n+2}}h_1(s)\,ds \cdot \int_{\eta_n}^{\eta_{n+2}}h_2(s)\,ds - 4\right) \label{eq51}\\
 &+ \int_{\widetilde\eta_{k+1}}^{\widetilde\eta_{k+3}}h_{1k}(s)\,ds \cdot \int_{\widetilde\eta_{k+1}}^{\widetilde\eta_{k+3}}h_{2k}(s)\,ds - 4\,.  \notag
\end{align}
Indeed, $\widetilde\Hh_{(k)}$ is constant on $[\eta_{k+2}, +\infty) = [\widetilde\eta_{k+2}, +\infty)$ and $\Hh = \widetilde\Hh_{(k)}$ on $[0, \eta_{n+2}]$, hence the terms with indexes $n \ge k+2$ in formula \eqref{eq2} for $\widetilde \Hh_{(k)}$ vanish, while the terms with indexes $n \le k$ coincide with the corresponding terms in \eqref{eq2} for the Hamiltonian $\Hh$. Since $\widetilde \Hh_{(k)} = \diag(a_{1k}, a_{2k})$ on $[\eta_{k+2}, + \infty)$, we have 
$$
\int_{\widetilde\eta_{k+1}}^{\widetilde\eta_{k+3}}h_{1k}\,ds \cdot \int_{\widetilde\eta_{k+1}}^{\widetilde\eta_{k+3}}h_{2k}\,ds 
= \prod_{j =1}^{2}\left(\int_{\eta_{k+1}}^{\eta_{k+2}}h_{j}\,ds + a_{jk}(\widetilde\eta_{k+3} - \widetilde\eta_{k+2})\right).
$$
A short calculation gives $\widetilde\eta_{k+3} - \widetilde\eta_{k+2} = 1/\sqrt{a_{1k} a_{2k}}$. Thus, we have 
$$
\int_{\widetilde\eta_{k+1}}^{\widetilde\eta_{k+3}}h_{1k}\,ds \cdot \int_{\widetilde\eta_{k+1}}^{\widetilde\eta_{k+3}}h_{2k}\,ds 
= \left(x_1 + \sqrt{\tfrac{a_{1k}}{a_{2k}}}\right)\left(x_2 + \sqrt{\tfrac{a_{2k}}{a_{1k}}}\right),
$$
where $x_j = \int_{\eta_{k+1}}^{\eta_{k+2}}h_{j}\,ds$ for $j =1,2$. Denoting $y_j = \int_{\eta_{k+2}}^{\eta_{k+3}}h_{j}\,ds$, $j =1,2$, we get 
\begin{equation}\label{eq52}
\left(x_1 + \sqrt{\tfrac{a_{1k}}{a_{2k}}}\right)\left(x_2 + \sqrt{\tfrac{a_{2k}}{a_{1k}}}\right) \le (x_1 + y_1)(x_2 + y_2) = \int_{\eta_{k+1}}^{\eta_{k+3}}h_{1}\,ds \cdot \int_{\eta_{k+1}}^{\eta_{k+3}}h_{2}\,ds,
\end{equation} 
for the following special choice of parameters $a_{1k}$ and $a_{2k}$: $a_{1k} = y_1^2$, $a_{2k} = 1$, where the  inequality in \eqref{eq52} follows from 
$y_1 y_2 \ge \bigl(\int_{\eta_{k+2}}^{\eta_{k+3}} \sqrt{h_1 h_2}ds\bigr)^2 = (\xi_{\Hh}(\eta_{k+3}) - \xi_{\Hh}(\eta_{k+2}))^2 = 1.$
Combining \eqref{eq51} and \eqref{eq52}, we see that $\widetilde \K(\widetilde\Hh_{(k)}) \le \widetilde \K(\Hh)$ for every $k$ and 
\begin{equation}\label{eq53}
\lim_{k \to \infty} \widetilde \K(\widetilde\Hh_{(k)}) = \widetilde\K(\Hh).
\end{equation}
By Lemma \ref{l0}, the spectral measure of the Hamiltonian $\widetilde \Hh_{(k)}$ belongs to $\sz$  for every $k$. From Lemma \ref{l6}, Lemma \ref{l7}, and \eqref{eq53} we obtain $\mu \in \sz$ and
$$
\K(\mu) \le \limsup_{k\to \infty} \K_{\widetilde\Hh_{(k)}}(0) \le c\limsup_{r\to \infty}\widetilde \K(\widetilde \Hh_{(k)})e^{c\widetilde \K(\widetilde \Hh_{(k)})} \le  c\widetilde \K(\Hh)e^{c\widetilde \K(\Hh)},
$$
with an absolute constant $c$. The theorem is proved. \qed 

\medskip

\section{Functions with summable fixed-scale Muckenhoupt characteristic}\label{s5}
In this section, we study functions from the class $A_2(\R_+, \ell^1)$ defined in Section \ref{s4} and prove Propositions \ref{p9}, \ref{p5}.  
\begin{Lem}\label{l51} Let $I=I^-\cup I^+$ be a splitting of an interval $I \subset \R$ into the union of two disjoint subintervals  $I^{\pm}$. Let $h \ge 0$ be a function on $I$ such that $h, 1/h \in L^1(I)$,  and let $\eta= \langle h\rangle_I  \langle 1/h\rangle_I-1$. Assume that $|I^-|/|I| \ge \tfrac{1}{5}$, then
\begin{equation}\label{cook131}
\left|\frac{\langle h\rangle_{I}}{\langle h\rangle_{I^-}}-1\right| \lesssim \sqrt{\eta(1+\eta)}, \quad 
\left|  \frac{\langle h\rangle_{I^-}}{ \langle h\rangle_{I}}-1  \right|\lesssim  \min(1, \sqrt{\eta}),
\end{equation}
and, moreover, 
\begin{equation}\label{sdp4}
\langle h\rangle_{I^-}\langle 1/h\rangle_{I^-}-1\lesssim  \eta. 
\end{equation}
\end{Lem}
\beginpf The number $\eta$ and all bounds are invariant with respect to multiplying $h$ with a positive constant, thus we can assume that $\langle h\rangle_{I}=1$. Next, put $\upsilon= |I^-|/|I|$, $a^{\pm}= \langle h\rangle_{I^\pm}$, $b^{\pm}= \langle h^{-1}\rangle_{I^\pm}$. We have
\begin{equation}\label{eq33}
\upsilon a^-+(1-\upsilon)a^+=1, \qquad \upsilon b^-+(1-\upsilon)b^+= \langle h^{-1}\rangle_I=1+\eta, \qquad a^{\pm}b^{\pm}\ge 1.
\end{equation}
Adding the first two estimates and using the bounds $ 1/a^{\pm}\le b^{\pm}$, one gets
 $\upsilon \left(a^- + 1/a^-\right) +(1-\upsilon)\left(a^+ + 1/a^+\right) \le 2 + \eta$.
Since $x+1/x \ge 2$ for all $x > 0$, this yields $\upsilon(a^- + 1/a^-) \le 2\upsilon +\eta$.
Dividing by $2v$, we get the inequality 
\begin{equation}\label{eq32}
\frac{1}{2}\left(a^- + \frac{1}{a^-}\right) \le 1 + \frac{\eta}{2\upsilon}.
\end{equation} 
It can be rewritten in the form $(1/a^- - 1)^2 \le  \eta/(\upsilon a^-)$. Since $\upsilon \in [\tfrac{1}{5}, 1]$ and $1/a^- \lesssim (1 + \eta)$ by \eqref{eq32}, this gives the first bound in \eqref{cook131}. 
To get the second bound in \eqref{cook131}, rewrite \eqref{eq32} in the form $(a^- - 1)^2 \le a^- \eta/\upsilon$ and use the fact that $\upsilon a^- \le 1$. Thus, 
\[
|a^--1|\le \frac{\sqrt{\eta}}{\upsilon},\, |a^--1|\le 1+\upsilon^{-1}\,,
\]
which implies the second inequality in \eqref{cook131}. Next, let us prove \eqref{sdp4}. Since $a^\pm + b^\pm \ge 2$, we get  $v(a^- + b^-) \le 2\upsilon + \eta$ by summing up the first two identities in \eqref{eq33}.
Hence $\sqrt{a^- b^-} \le 1 + \eta/(2\upsilon)$ and $a^-b^-\le 1+\eta/\upsilon+\eta^2/(4\upsilon^2)$.
This gives the inequality 
$ \langle h\rangle_{I^-}\langle 1/h\rangle_{I^-}-1\lesssim  \eta$ in the case where $\eta \le \upsilon$. For $\eta \ge \upsilon$ we can use \eqref{eq33} to get $a^- \le 1/\upsilon \le 5$ and $b^{-} \le 5 (1+\eta)$. This gives 
$\langle h\rangle_{I^-}\langle 1/h\rangle_{I^-}-1 \le 25(1+\eta)-1\lesssim \eta$ since $\eta\ge 1/5$.  \qed

\medskip

\noindent {\bf Proof of Proposition \ref{p9}.} Apply Lemma \ref{l51} to the function $h$ and the intervals $I = I_{n, \alpha_n}$, $I_- = [n, n+2]$, $n \ge 0$. Since $\{\alpha_n\} \subset [3, 4]$, this will give the estimate $[h]_{2, \ell^1} \le c[h, \alpha]$ with an absolute constant $c$. \qed   

\medskip
\begin{Lem}\label{l52} For $h\in A_2(\mathbb{R}_+,\ell^1)$, define
$Q_n = \langle h\rangle_{I_{n,2}}\langle h^{-1}\rangle_{I_{n,2}}-1$ and $f_n = \langle h\rangle_{I_{n,1}}$. Then, 
\begin{equation}
 (1+Q_n)^{-1} \lesssim \frac{f_{n+1}}{f_n}\lesssim 1+Q_n, \label{cookk9}
\end{equation}
\begin{equation}\label{cook9}
 \left|\frac{f_{n+1}}{f_n} - 1 \right|\le c\sqrt{Q_n},\;\; \mbox{if} \;\; Q_n \le 1\,.
\end{equation}
 Moreover, we have $\|\widetilde h+\widetilde h^{-1}-2\|_1\lesssim [h]_{2,\ell^1}=\sum_{n=0}^\infty Q_n$ for the function $\widetilde h$ defined by 
\begin{equation}\label{sdl1}
\widetilde h(x)=h(x)/\langle h\rangle_{I_{n,1}}, \quad x\in I_{n,1}, \quad n\in \mathbb{Z}^+.
\end{equation}
\end{Lem}
\beginpf 
Represent $f_{n+1}/f_{n}$ in the form
\begin{equation}\label{eq35}
\frac{f_{n+1}}{f_{n}} 
= \frac{\langle h \rangle_{I_{n+1,1}}}{\langle h \rangle_{I_{n, 2}}}
\frac{\langle h \rangle_{I_{n, 2}}}{\langle h \rangle_{I_{n, 1}}}\,.
\end{equation}
We write
\begin{equation}\label{sdl4}
\frac 12\le    \frac{\langle h \rangle_{I_{n,2}}}{\langle h \rangle_{I_{n, 1}}}
\le  1+ c\sqrt{Q_n(Q_n + 1)}\lesssim 1+Q_n,
\end{equation}
where the first inequality is immediate and the second one follows from the  first estimate in \eqref{cook131}.  Similarly, we get
$$
\frac 12\le    \frac{\langle h \rangle_{I_{n,2}}}{\langle h \rangle_{I_{n+1, 1}}}
\le  1+ c\sqrt{Q_n(Q_n + 1)}\lesssim 1+Q_n
$$
and
\begin{equation}\label{sdl3}
( 1+ Q_n)^{-1}         \lesssim    \frac{\langle h \rangle_{I_{n+1,1}}}{\langle h \rangle_{I_{n, 2}}}
\le 2\,.
\end{equation}
It is now sufficient to multiply \eqref{sdl3} with \eqref{sdl4} and substitute into \eqref{eq35} to get \eqref{cookk9}.
Take $n \ge 0$ such that $Q_n \le 1$. By Lemma \ref{l51}, we have 
\begin{equation}\label{eq34}
\left|\frac{\langle h \rangle_{I_{n, 2}}}{\langle h \rangle_{I_{n,1}}} - 1 \right| \lesssim \sqrt{Q_n},\quad
\left|\frac{\langle h \rangle_{I_{n+1,1}}}{\langle h \rangle_{I_{n, 2}}} - 1 \right| \lesssim \sqrt{Q_n} \,.
\end{equation}
\medskip
Substituting these bounds into \eqref{eq35} gives \eqref{cook9}.
Finally, observe that for every $n \ge 0$ we have
$\langle h\rangle_{I_{n,1}}\langle h^{-1}\rangle_{I_{n,1}}-1 \lesssim Q_n$ by \eqref{sdp4}.
Using the identity
\[
\sum_{n=0}^\infty \|\widetilde h+\widetilde h^{-1}-2\|_{L^1(I_{n,1})}=2\sum_{n=0}^\infty\left(\langle h\rangle_{I_{n,1}}\langle h^{-1}\rangle_{I_{n,1}}-1\right),
\]
we complete the proof of the lemma. \qed

\medskip

\noindent {\bf Remark.} Notice that \eqref{cookk9} and \eqref{cook9} imply 
\begin{equation}\label{cook42}
|\log (f_{n+1}/f_n)|\lesssim 
\begin{cases}
\sqrt Q_n, & Q_n<2\,,\\
\log Q_n, & Q_n>2.
\end{cases}
\end{equation}

\medskip

\noindent {\bf Proof of Proposition \ref{p5}.} Define $\widetilde h$ as in  \eqref{sdl1} and consider the function 
$f_1 = (\widetilde h-1)\chi_{\frac{1}{2}<\widetilde h<\frac{3}{2}}$. For shorthand,  denote $P=[h]_{2,\ell^1}=\sum_{n=0}^\infty Q_n$ where $Q_n$ is defined in the previous lemma. Since the function $\widetilde h + \widetilde h^{-1} -2 \in L^1(\R_+)$, we have $f_1 \in L^2(\R_+)$ and $\|f_1\|_{2}^{2} \lesssim P$. Indeed, this follows from the fact that $x + x^{-1} -2\sim (x-1)^2$ for $x \in [\frac{1}{2}, \frac{3}{2}]$ and the estimate $\|\widetilde h + \widetilde h^{-1} -2\|_{1} \lesssim P$ in Lemma~\ref{l52}. Similarly, the function $f_2 = (\widetilde h-1)\chi_{|\widetilde h-1| \ge \frac{1}{2}}$ belongs to $L^1(\R_+)$ and $\|f_2\|_{1} \lesssim P$. Thus, we see that $\widetilde h$ can be represented in the form $\widetilde h = f_0 + f_1 + f_2$, where $f_0=1$,
$f_1 \in L^2(\R_+)$, $f_2 \in L^1(\R_+)$, and $\|f_1\|_{2}^{2} + \|f_2\|_{1} \lesssim P$. Function $\widetilde h^{-1}$ admits similar  representation $\widetilde h^{-1} = \widehat f_0 + \widehat f_1 + \widehat f_2$, where $\widehat f_0=1$, $\widehat f_1 = - f_1$ and $\widehat f_2 \in L^1(\R_+)$ is such that $\|\widehat f_2\|_{1} \lesssim P$. Notice that we have got $\widehat f_1=-f_1$ from 
\[
\frac{\chi_{|\widetilde h-1|<1/2}}{\widetilde h}=\frac{\chi_{|\widetilde h-1|<1/2}}{1+f_1}=\chi_{|\widetilde h-1|<1/2}(1-f_1+O(f_1^2))
\]
and  $\widehat f_2\in L^1(\mathbb{R}_+)$ because  $\widehat f_2=\chi_{|\widetilde h-1|<1/2}O(f_1^2)+\chi_{|\widetilde h-1|>1/2}(\widetilde h^{-1}-1)\in L^1(\mathbb{R}_+)$.

Let $g_0$ be the function on $\R_+$ such that $g_0=\log f_n$ on each $I_{n,1}$,  then $h = e^{g_0} \widetilde h$ on $\R_+$.  Define also  the function $g: x \mapsto g_0(x)-g_0(0)$ on $\R_+$. Then, for $\kappa$ and $\kappa_d$ from Proposition \ref{p141}, we have
\begin{align*}
\kappa &=\sum_{0 \le k, j \le 2}a_{ kj},  &a_{ kj}: x \mapsto \int_{x}^{\infty}\widehat f_{k}(x) f_j(\xi) e^{g(\xi) - g(x) + x - \xi}\,d\xi,
\\
\kappa_d &=\sum_{0 \le k, j \le 2}a_{d,kj},  &a_{d,kj}: x \mapsto \int_{x}^{\infty}f_{k}(x)\widehat f_j(\xi) e^{g(x) - g(\xi) + x - \xi}\,d\xi.\end{align*}
We will need some estimates for the function $g$. Let   $Q_j$, $f_j$ be defined as in Lemma \ref{l52}  and let $v_n= \log\bigl(f_{n}/f_{n-1}\bigr), n\in \mathbb{N}$, $v_0=0$. Observe that $g(x) = \sum_{n=0}^{[x]}v_n$ on $\R_+$ by construction. Here, as usual, $[x]$ stands for the integer part of a number $x \in \R_+$. 
We can estimate
\begin{equation}\label{cook65}
\|\{v_n\}\|_2^2=\sum_{n:\, Q_{n-1}<2}v_n^2+\sum_{n:\, Q_{n-1}>2}v_n^2\lesssim \sum_{n:\, Q_{n-1}<2}Q_n+\sum_{n:\, Q_{n-1}>2}\log^2 Q_n\lesssim P\,,
\end{equation}
where we used \eqref{cook42} and the trivial bound: $\log^2 Q\lesssim Q$ which holds for all $Q>2$. Bound \eqref{cook42} also yields
\begin{equation}\label{cook66}
\|\{v_n\}\|_\infty\lesssim \log(2+P)\,.
\end{equation}
For $x< y$, we can apply \eqref{cook42} to write
\begin{align}\label{cook76}
|g(x) - g(y)| \le \left|\sum_{j=[x]}^{[y]}v_j\right|
&\le \sum_{j=[x],\,Q_{j-1}<2}^{[y]}|v_j| + \sum_{j=[x],\,Q_{j-1}\ge 2}^{[y]}|v_j|\\     \nonumber 
&\lesssim \sum_{j=[x],\,Q_{j-1}<2}^{[y]}\sqrt{|Q_{j-1}|} + \sum_{j=[x],\,Q_{j-1}\ge 2}^{[y]}\log Q_{j-1}\nonumber\\
&\lesssim \Bigl((|x-y|+ 1)\sum_{j \ge 0} Q_j\Bigr)^{1/2} + \sum_{j\ge 0}Q_j\nonumber\\
&\lesssim \sqrt{(|x-y|+ 1)P} + P\,.\nonumber
\end{align}
It follows that there is an absolute constant $C$ such that for all $x,y \in \R_+$ we have
\begin{align}
|g(x) - g(y)| &\le \tfrac{1}{2}|x-y| + C(1+P). \label{eq36}
\end{align}
Now, for indexes $k$, $j$ such that $k+j\ge 2$, we can use \eqref{eq36} and the Young  inequality for convolutions to estimate
\begin{align*}
\|a_{d,kj}\|_{1}
&\lesssim e^{CP}\int_{0}^{\infty}\!\!\!\int_{0}^{\infty}\!\! |f_{k}(x)| \chi_{\R_+}(\xi-x)e^{-(\xi - x)/2} |\widehat f_j(\xi)|\,d\xi\,dx\\
&\lesssim e^{CP}\|f_{k}\|_{p_k} \cdot \|\chi_{\R_+}e^{-x}\|_{r_{k,j}} \cdot \|\widehat f_j\|_{p_j}
\lesssim P e^{CP}, 
\end{align*}
where $p_0 = +\infty$, $p_1 = 2$, $p_2 = 1$, and the parameter $r_{k,j}$ is chosen so that 
$\frac{1}{p_k} + \frac{1}{r_{k,j}} + \frac{1}{p_j}= 2$. The estimate on $a_{kj}$ for $k+j\ge 2$ is similar.
To prove that $\kappa + \kappa_d -2 \in L^1(\R_+)$, it remains to estimate the $L^1(\R^+)$--norms of functions
\begin{align*}
a_{00} + a_{d, 00} - 2&= 2\int_{x}^{\infty}e^{x - \xi} \left(\cosh G(x,\xi) - 1\right)\,d\xi, \\ 
a_{01} + a_{d, 01} &= 2\int_{x}^{\infty} \widehat f_1(\xi) e^{x - \xi} \sinh G(x,\xi) \,d\xi, \\  
a_{10} + a_{d, 10} &= 2\int_{x}^{\infty} f_{1}(x) e^{x - \xi} \sinh G(x,\xi)\,d\xi,
\end{align*}
where $G(x,\xi) = g(x) - g(\xi)$. Let us define the function $\widetilde g$  on $[-1,\infty)$ to be continuous, linear on $I_{j,1}$ for each $j\ge -1$, and so that $\widetilde g(-1)=0$, $\widetilde g(j) = \sum_{n=0}^{j}|v_n|$ for $j \ge 0$. Clearly, $\widetilde g$ is non-decreasing on $[-1,\infty)$. Put $\widetilde G(x,\xi) = \widetilde g(\xi+1) - \widetilde g(x-1)$ for every $0<x<\xi$. Then $|G(x,\xi)| \le \widetilde G(x,\xi)$ and so $\cosh G(x,\xi)\le \cosh \widetilde G(x,\xi)$. By construction and \eqref{cook65}, we have 
\begin{equation}\label{cook88}
\|\widetilde g'\|_{2}^{2} \lesssim \sum_{n \ge 0}|v_n|^2 \lesssim P.
\end{equation}
The bound \eqref{cook65} also implies 
\begin{equation}\label{cook54}
\|\widetilde G(x,x)\|_2^2 \lesssim \|\{v_n\}\|_2^2\lesssim P\,.
\end{equation}
 The estimate \eqref{cook66} gives 
\begin{equation}\label{eq55}
\|\widetilde G(x,x)\|_{\infty} \lesssim \sup_{n \ge 0}|v_n| \lesssim \log(2+P)
\end{equation}
and argument given in \eqref{cook76} yields
\begin{equation}\label{cook53}
\qquad \widetilde G(x,\xi)    \lesssim \sqrt{(|x-\xi|+ 1)P} + P, \quad   \widetilde G(x,\xi)  \le  \tfrac{1}{2}|x-\xi| + C(1+P)
\end{equation} 
for all $x<\xi$. Integrate by parts to get
\begin{align*}
\|a_{00} + a_{d, 00} - 2\|_{1} 
&\le 2\int_{0}^{\infty}\!\!\int_{x}^{\infty}e^{x-\xi}(\cosh\widetilde G(x, \xi) - 1)\,d\xi\,dx\\
&\le 2\int_{0}^{\infty}\!\!\int_{x}^{\infty}\widetilde g'(\xi+1)e^{x-\xi}\sinh\widetilde G(x, \xi)\,d\xi\,dx + 2R_1,  
\end{align*}
where $R_1 = \int_{0}^{\infty}(\cosh\widetilde G(x,x) - 1)\,dx$. Using the inequality $\cosh t - 1 \lesssim t^2e^{|t|}$, we obtain 
$R_1 \le \|\widetilde G(x,x)\|_2^2 \exp(\|\widetilde G(x,x)\|_\infty) \lesssim Pe^{CP}$ by \eqref{cook54} and \eqref{eq55}. To estimate the double integral, let us change the order of integration and integrate by parts once again: 
\begin{equation}\label{cook08}
\int_{0}^{\infty} \widetilde g'(\xi+1) \int_{0}^{\xi} e^{x-\xi}\sinh\widetilde G(x, \xi)\,dxd\xi
=\int_{0}^{\infty}\!\!\widetilde g'(\xi+1) \int_{0}^{\xi}\widetilde g'(x-1) e^{x-\xi}\cosh\widetilde G(x, \xi)\,dx\,d\xi + R_2,
\end{equation}
where $R_2 = \int_{0}^{\infty}\widetilde g'(\xi+1) (\sinh \widetilde G(\xi, \xi) - e^{-\xi} \sinh \widetilde G(0, \xi))\,d\xi\le  \int_{0}^{\infty}\widetilde g'(\xi+1) \sinh \widetilde G(\xi, \xi) d\xi$ because $\widetilde g'\ge 0$.
Let us estimate the integral first using the second  bound in \eqref{cook53}
\begin{align*}
\int_{0}^{\infty}\!\!\widetilde g'(\xi+1) \int_{0}^{\xi}\widetilde g'(x-1) e^{x-\xi}\cosh\widetilde G(x, \xi)\,dx\,d\xi 
&\lesssim e^{CP} \int_{0}^{\infty}\!\!\widetilde g'(\xi+1) \int_{0}^{\xi}\widetilde g'(x-1) e^{(x-\xi)/2}\,dx\,d\xi \\ 
&\lesssim e^{CP}\|\widetilde g'\|_2^2\lesssim Pe^{CP},
\end{align*}
as follows from Young's inequality for convolution and \eqref{cook88}. We are left with estimating $R_2$. Using  inequality  $|\sinh t| \le |t| e^{|t|}$ we obtain
$$
 \int_{0}^{\infty}\widetilde g'(\xi+1)   \sinh \widetilde G(\xi, \xi)       d\xi
 \le \|\widetilde g'(\xi+1)\|_{2}\cdot\|\widetilde G(\xi,\xi)\|_{2}\exp(\|\widetilde G(\xi,\xi)\|_\infty) \lesssim Pe^{CP}\,.$$
Collecting the bounds, we get $\|a_{00} + a_{d, 00} - 2\|_1\lesssim Pe^{CP}$. It remains to bound the $L^1(\R_+)$--norms of $a_{01} + a_{d, 01}$ and $a_{10} + a_{d, 10}$. First, we write
$$
\|a_{01} + a_{d, 01}\|_{1} \le 2\int_{0}^{\infty}|\widehat f_1(\xi)|\!\!\int_{0}^{\xi}e^{x-\xi}\sinh\widetilde G(x,\xi)\,dx\,d\xi\lesssim Pe^{CP}
$$
since the integral has the form similar to the left hand side in \eqref{cook08} and  the estimates for \eqref{cook08} can be repeated. Finally,
\begin{align*}
\|a_{10} + a_{d, 10}\|_{1} 
\le &2\int_{0}^{\infty}\!\!\int_{x}^{\infty} |f_{1}(x)| e^{x - \xi} \sinh \widetilde G(x,\xi)\,d\xi\,dx\\
\le &2\int_{0}^{\infty}|f_{1}(x)|\sinh \widetilde G(x,x)\,dx \\
    &+ 2\int_{0}^{\infty}\!\!\int_{x}^{\infty} |f_{1}(x)| \widetilde g'(\xi+1) e^{x - \xi} \cosh \widetilde G(x,\xi)\,d\xi\,dx,
\end{align*}
where the first term can be estimated similarly to $R_2$, while the second one is dominated by 
\mbox{$Ce^{CP}\|f_1\|_{2}  \cdot \|\widetilde g'(t-1)\|_{2}$} $\lesssim Pe^{CP}$. Thus, we see that $\kappa + \kappa_d - 2$ belongs to $L^1(\R_+)$ and $[h]_{int} \lesssim  P e^{cP}$ with an absolute constant~$c$. \qed

\medskip

\section{Krein strings and proof of Theorem \ref{t2}}\label{s6}
In this section, we introduce the spectral measure for Krein string and show how Theorem \ref{t1} and some results obtained in \cite{KWW} imply Theorem \ref{t2}. Let $0 < L \le \infty$. Recall that $M$ and $L$ form $[M,L]$ pair if \eqref{sdsing} holds, i.e.,
$L+\lim_{t\to L}M(t)=\infty\,\quad {\rm and}\,\quad \lim_{t\to L}M(t)>0$.
Define the Lebesgue--Stieltjes measure $\mf$ by $\mf[0,t]=M(t)$.
Next, define the increasing function  $N: t \mapsto t+M(t)$ on $[0, L)$ and let $\nf$ denote the corresponding measure, $\nf[0,t]=N(t)$ for $t \ge 0$. Define also the function $N^{(-1)}$ on $\R_+$ by $N^{(-1)}: y \mapsto  \inf \{ t\ge 0: N(t) \ge y \}$. The set under the last infimum is non-empty for every $y \ge 0$ because of the assumptions we made on $M$ and $L$. Using the fact that $N$ is strictly increasing, one can show that $N^{(-1)}$ is continuous on $\R_+$, and we have $N^{(-1)}(N(t)) = t$ for every $t \in [0, L)$. Let $M'$ be the density of the absolutely continuous part of $\mf$, so that $\mf=M'(t)\,dt+\mf_s$. Denote by $E_s$ the support of the singular part $\mf_s$ of the measure~$\mf$. Define two functions on $\R_+$, 
\begin{equation}
h_1(x)= 
\begin{cases}
0, & \mbox{if } N^{(-1)}(x)\in E_s,\\
\frac{1}{1+M'(N^{(-1)}(x))}, & \mbox{otherwise},
\end{cases}
\label{sd1}
\end{equation}
and
\begin{equation}
h_2(x)= 
\begin{cases}
1, & \mbox{if } N^{(-1)}(x)\in E_s,\\
\frac{M'(N^{(-1)}(x))}{1+M'(N^{(-1)}(x))}, & \mbox{otherwise}.
\end{cases}
\label{sd2}
\end{equation}
The proof of Lemma~\ref{l20} below shows that functions $h_{1}$, $h_2$ defined by different representatives of the function $M'$ differ on a set of zero Lebesgue measure. Notice that $h_{1}$, $h_2$ are non-negative Lebesgue measurable functions and we have $h_1(x)+h_2(x)=1$ for  all $x \in \R_+$. 
We are going to prove the following result from \cite{KWW}, pp. 1527--1528.
\begin{Lem}\label{l20}
Formulas \eqref{sd1}, \eqref{sd2} establish the bijection $[M ,L] \mapsto \diag(h_1,h_2)$ between $[M, L]$ pairs and nontrivial diagonal Hamiltonians $\Hh = \diag(h_1, h_2)$ with unit trace almost everywhere on~$\R_+$.
\end{Lem}
\beginpf 
Fix any pair $[M, L]$ and consider the corresponding function $N^{(-1)}$ and the measure $\nf$. For  every function $f \in L^1_{\rm loc}(\R_+,\nf)$ we have $f(N^{(-1)}(x))\in L^1_{\rm loc}(\R_+)$, and, moreover,
\begin{equation}\label{kor}
\int_{[0, L)} f(t)\,d\nf(t)=\int_{\R_+} f(N^{(-1)}(x))\,dx,
\end{equation} 
if $f$ is compactly supported in $[0, L)$. This result is known as the change of variables in the Lebesgue--Stieltjes integral (see, e.g., Exercise 5 in Section III.13 of\cite{Dunford58}) but we give its proof for completeness. Without loss of generality we can assume that $f\ge 0$. Then (see, e.g., \cite{KWW}, 
Proposition 6.24), we have
\[\int_{[0, L)} f(t)\,d\nf(t)=\int_{\R_+} \Lambda_1(\lambda)\,d\lambda, \qquad
\int_{[0,L)} f(N^{(-1)}(x))\,dx=\int_{\R_+} \Lambda_2(\lambda)\,d\lambda,
\] 
where $\Lambda_1(\lambda) = \nf\{t: f(t)>\lambda\}$ and  $\Lambda_2(\lambda) = \left|\{x: f(N^{(-1)}(x))>\lambda\}\right|$. For all $0\le a<b$ we have 
\begin{equation}\label{sdnep}
\nf((a,b))=N(b-)-N(a)=|(N(a),N(b-))|\,,
\end{equation}
where $N(b-)$ denotes the left limit of $N$ at the point $b$. In fact, $(N(a), N(b-))$ is preimage of $(a,b)$ under the continuous map $N^{(-1)}$.  Thus, the preimage under $N^{(-1)}$ of any open cover $\cup (a_j, b_j)$ for  $\nf$-measurable set $E$ will be an open cover for the set $\{x: N^{(-1)}(x)\in E\}$. Conversely, every open cover $\cup_j (c_j,d_j)$ for $\{x: N^{(-1)}(x)\in E\}$ is the preimage of some open cover for $E$. Indeed, for each $j$ we get $(c_j,d_j)=(N(a_j),N(b_j))$, where $a_j$ and $b_j$ are points of continuity for $N$ (to see this, note that the preimage of $\nf$'s atom under $N^{(-1)}$ is a closed segment). For every regular measure $\nu$ we have
\begin{equation}\label{sdnep1}
\nu(E)=\inf\Bigl\{\sum_{j}\nu(I_j),\;  E\subset \cup_j I_j, \, \{I_j\} \mbox{ are disjoint  open intervals}\Bigr\},
\end{equation}
see, e.g., Lemma 1.17 in \cite{folland}. From \eqref{sdnep} and \eqref{sdnep1} we now get $\Lambda_1(\lambda)=\Lambda_2(\lambda)$ and, consequently, relation \eqref{kor} follows. Next, take a number $y \ge 0$. Since $h_1(x) = 0$ for all $x$ such that $N^{(-1)}(x) \in E_s$, we have 
$$
\chi_{[0,y]}(x)h_1(x) = f_y(N^{(-1)}(x)), \quad x \in [0,L), 
$$  
where $f_y: t \mapsto \frac{\chi_{[0, N^{(-1)}(y)] \setminus E_s}(t)}{1+M'(t)}$ is the compactly supported function from $L^{1}([0, L), \nf)$. Applying formula \eqref{kor} to the function $f_y$, we get
\begin{equation}\label{e1} 
\int_0^y h_1(x)dx 
=\int_{[0,L)}\frac{\chi_{[0, N^{(-1)}(y)] \setminus E_s}(t)}{1 + M'(t)}\,d\nf(t) 
= \int_{[0, N^{(-1)}(y)] \setminus E_s}dt = N^{(-1)}(y),
\end{equation}
where we used the fact that the singular part of $\nf$ is supported on $E_s$ and the absolutely continuous part of $\nf$ has density $M' +1$ with respect to the Lebesgue measure on $[0,L)$. If $y$ is a point of growth for the function $N^{(-1)}$ (that is, there is no open interval $I$ containing $y$ such that $N^{(-1)}$ is constant on $I$), we have $\chi_{[0,y]}(x) = \chi_{[0,N^{(-1)}(y)]}(N^{(-1)}(x))$ for all $x \ge 0$, hence we can apply \eqref{kor} to get
\begin{equation} \label{eq421}
\int_0^y h_2(x)dx=\int_{[0, N^{(-1)}(y)]\setminus E_s}\frac{M'(t)(1+M'(t))}{1+M'(t)}dt+\int_{[0, N^{(-1)}(y)] \cap E_s}\!\!d\mf_s = \mf[0,N^{(-1)}(y)].  
\end{equation}
From here we see that $h_1$, $h_2$  define $M$, $L$ uniquely, in particular, these functions, as elements of $L^1_{\rm loc}(\mathbb{R}_+)$, do not depend on the choice of the representative of $M'$. Moreover, we cannot have $h_1=0$ or $h_2=0$ almost everywhere on $\R_+$ for any $M$, $L$ satisfying \eqref{sdsing}. Hence, $[M, L] \mapsto \diag(h_1, h_2)$ is the injective mapping from a set of pairs $[M,L]$ to nontrivial diagonal Hamiltonians with unit trace. Now take a nontrivial Hamiltonian $\diag(h_1, h_2)$ with unit trace almost everywhere on~$\R_+$, and consider the function 
$$
\Psi: y \mapsto \int_0^y h_1(x)\,dx.  
$$
Put $L = \sup_{y\ge 0} \Psi(y)$. Note that $|\Psi(y_1) - \Psi(y_2)| \le |y_1 - y_2|$ for all $y_1$, $y_2$ in $\R_+$, hence there exists a measure $\mf$ on $[0, L)$ such that $\Psi(y) =  \inf \{ x\ge 0: x + M(x) \ge y\}$ for every $y \ge 0$, where $M(x) = \mf[0, x]$. Using \eqref{e1} and \eqref{eq421}, it is easy to check that formulas \eqref{sd1}, \eqref{sd2} for $[M, L]$ generate the singular Hamiltonian $\Hh = \diag(h_1, h_2)$ and it is nontrivial. The lemma is proved. \qed

\medskip

For any pair $[M,L]$, one can define the  Krein string as the differential operator \cite{KK68, DymMcKean}. 
In \cite{KWW}, the authors considered two functions $\phi(x,z)$ and $\psi(x,z)$ that satisfy  
\[	
\phi(x,z)=1-z\int_{[0,x]}(x-s)\phi(s,z)\,d\mf (s), \quad x\in [0,L)\,,
\]
\[
\psi(x,z)=x-z\int_{[0,x]}(x-s)\psi(s,z)\,d\mf (s), \quad x\in [0,L)\,.
\]
These functions are uniquely determined by the string $[M,L]$ and they define the principal Weyl-Titchmarsh function $q$ of $[M,L]$ by
\[
q(z)=\lim_{x\to L} \frac{\psi(x,z)}{\phi(x,z)}, \quad z\in \mathbb{C}\backslash [0,\infty),
\]
see formula (2.21) in \cite{KWW}. This function $q$ has the unique integral representation 
\[
q(z)=b+\int_{\R_+} \frac{d\sigma(x)}{x-z}\,,
\]
where $b\ge 0$ and $\sigma$, the spectral measure of the string $[M,L]$, is a measure on $\R_+ = [0,+\infty)$ satisfying condition
\[
\int_{\R_+}\frac{d\sigma(x)}{1+x}<\infty\,.
\]
The authors of \cite{KWW} established, among other things, connection between $q$ and the Weyl-Titchmarsh function of a canonical system.
It is worth to mention that the definition of the Weyl-Titchmarsh function $m$ we used in \eqref{eq9} was taken from \cite{Romanov}. The authors of \cite{HSW}, \cite{KWW}  deal with the canonical system written differently, i.e., they write the Cauchy problem
$$
W'(t,z)J = z W(t,z)\Hh(t), \qquad W(0,z) = \idm, \qquad t \in \R_+, \quad z \in \C,
$$ 
and define the Weyl-Titchmarsh function $Q^+$ for $z \in \C\setminus\R$ by
\begin{equation}\label{eq70}
Q^+(z) = \lim_{t \to +\infty}\frac{w_{11}(t,z)\tilde\omega + w_{12}(t,z)}{w_{21}(t,z)\tilde\omega + w_{22}(t,z)}, 
\qquad
W(t,z) = \begin{pmatrix}w_{11}(t,z)&w_{12}(t,z)\\w_{21}(t,z)&w_{22}(t,z)\end{pmatrix}.
\end{equation} 
It is not difficult to see that $W(t,z) = M(t, -z)^\top$ for the solution $M$ of \eqref{eq1}. If we let $\sigma_1 = \left(\begin{smallmatrix}0&1\\1&0\end{smallmatrix}\right)$ and denote by $M_{\sigma_1}$ the solution of Cauchy problem $JM'_{\sigma_1} = z \Hh_{\sigma_1} M_{\sigma_1}$, $M_{\sigma_1}(0,z) = \idm$ for the dual  Hamiltonian $\Hh^d=\Hh_{\sigma_1} = \sigma_1\Hh\sigma_1$, then the function $m_{\sigma_1}$ from formula \eqref{eq9} for $\Hh_{\sigma_1}$ will coincide with the function $Q^+$ in \eqref{eq70} for $\Hh$ and $\tilde\omega = 1/\omega$. Indeed, we have 
\begin{equation}\label{sdfe3}
M_{\sigma_1}(t,z) = \sigma_1 M(t, -z)\sigma_1 = \sigma_1 W(t,z)^{\top}\sigma_1 = \begin{pmatrix}w_{22}(t,z)&w_{12}(t,z)\\w_{21}(t,z)&w_{11}(t,z)\end{pmatrix}.
\end{equation}
We will need the following lemma from \cite{KWW}.
\begin{Lem}
Suppose $[M,L]\mapsto \diag(h_1,h_2)$ is the bijection given by \eqref{sd1} and \eqref{sd2}, $q$ is the Weyl-Titchmarsh function for the string given by $[M,L]$, and $m,m_{\sigma_1}$ are the Weyl-Titchmarsh functions for $\diag(h_1,h_2)$ and  $\diag(h_2,h_1)$, respectively. Then, we have 
\begin{equation}\label{sdk1}
zq(z^2)=m_{\sigma_1}(z)=-m^{-1}(z), \quad z\in \mathbb{C}^+\,.
\end{equation}
\end{Lem}
\beginpf In \cite{KWW}, formula (4.20), it is proved that 
\begin{equation}
Q^+(z)=zq(z^2), \quad z\in \mathbb{C}^+\,,
\end{equation}
where $Q^+$ is defined in \eqref{eq70} and $\Hh$ is obtained from $[M,L]$ by bijection discussed in Lemma \ref{l20}. On the other hand, $Q^+(z)=m_{\sigma_1}(z)=m^{-1}(-z)=-m^{-1}(z)$, where the first equality follows from discussion right before formula \eqref{sdfe3}, the second one follows from \eqref{sdfe3} and \eqref{eq9}, and the last one is the corollary of the spectral measure of  $\diag(h_1,h_2)$ being even. \qed

\medskip

\noindent{\bf Proof of Theorem \ref{t2}.} Let $[M, L]$ be a  string with Weyl-Titchmarsh function $q$ and the spectral measure $\sigma$. Using Lemma \ref{l20}, define the Hamiltonians $\Hh$ and $\Hh^d=\Hh_{\sigma_1} = \sigma_1 \Hh \sigma_{1}$ on $\R_+$. Let $m_{\sigma_1}$, $\mu_{\sigma_1} = w_{\sigma_1}\,dx + \mu_{\sigma_1,s}$ be the Weyl-Titchmarsh function and the spectral measure of $\Hh^d$. Recall that $\sigma=v\,dx+\sigma_s$ for  spectral measure of the string. In \eqref{sdk1}, taking the nontangential limits  of $\Im(m_{\sigma_1}(z))$ and $\Im (zq(z^2))$   as $z \to x$, we get $w_{\sigma_1}(x)$ and $xv(x^2)$ for almost all $x \in \R_+$, respectively. Thus, $w_{\sigma_1}(x)=xv(x^2)$ for almost every $x\ge 0$, and, since $\mu_{\sigma_1}$ is even by Lemma \ref{l11}, we get
\[
\int_{\mathbb{R}} \frac{\log w_{\sigma_1}(x)}{1+x^2}\,dx= 2\int_0^\infty \frac{\log x}{1+x^2}\,dx+2\int_0^\infty \frac{\log v(x^2)}{x^2+1}\,dx=\int_0^\infty \frac{\log v(x)}{\sqrt x(x+1)}\,dx,
\]
where we used the fact that $\int_0^\infty \frac{\log x}{1+x^2}\,dx=\int_{-\infty}^{+\infty} \frac{y}{e^y+e^{-y}}\,dy=0$.
This implies that $\int_0^\infty \frac{\log v(x)}{\sqrt x(x+1)}\,dx$ is finite if and only if $\mu_{\sigma_1} \in \sz$. On the other hand, formula \eqref{kor} and the defintion of $h_1$, $h_2$ imply
$$
\int_{0}^{y}\sqrt{h_1(x)h_2(x)}\,dx = \int_{[0,N^{(-1)}(y)]\setminus E_s}\frac{\sqrt{M'(t)}}{1+M'(t)}\,d\nf(t) 
= \int_{0}^{N^{(-1)}(y)}\sqrt{M'(t)}\,dt 
$$
if $y$ is a point of growth of the function $N^{(-1)}$. For every $n \ge 1$ the points $\{\eta_n\}$ defined in \eqref{eq71} are the points of growth for $N^{(-1)}$. Indeed, this is clear from the formula \eqref{e1} that was proved for all $y \ge 0$. Hence we have $t_n = N^{(-1)}(\eta_n)$ for all $n \ge 0$. It follows that
$$
t_{n+2} - {t_n} = N^{(-1)}(\eta_{n+2}) - N^{(-1)}(\eta_{n}) = \int_{\eta_n}^{\eta_{n+2}}h_1(x)\,dx,
$$
where we used \eqref{e1} again. We also have
$$
M(t_{n+2})-M(t_n) = \mf(t_n, t_{n+2}] = \mf(N^{(-1)}(\eta_n), N^{(-1)}(\eta_{n+2})] = \int_{\eta_{n}}^{\eta_{n+2}}h_2(x)\,dx, 	
$$
by the definition of $M$ and \eqref{eq421}. Thus, $\widetilde\K[M, L] = \widetilde\K(\Hh)=\widetilde\K(\Hh_{\sigma_1})$ and $\sqrt{\det\Hh} \in L^1(\R_+)$ if and only if $\sqrt{M'} \in L^1(\R_+)$. Now the result follows from Theorem \ref{t1}. \qed

\medskip

\noindent {\bf Remark.} If $[M,L]\mapsto \diag(h_1,h_2)$, then the string $[M_d,L_d]$ for which  $[M_d,L_d]\mapsto \diag(h_2,h_1)$ is called the dual string. One can easily see that $\widetilde\K[M, L]=\widetilde\K[M_d, L_d]$ so the logarithmic integral for the string converges if and only it converges for the dual string.\medskip

We give two applications of Theorem \ref{t2}.

\begin{Prop} Suppose that the mass distribution $M$ of a string $[M, \infty]$ satisfies $M' = 1$ almost everywhere on $\R_+$. Let
$\mf_s$ be the singular measure on $\R_+$ such that $M(t) = t + \mf_s[0,t]$ for all $t \ge 0$. 
Then we have 
$$
\int_0^\infty \frac{\log v(x)}{\sqrt x(x+1)}\,dx > -\infty
$$ 
for the  spectral measure $\sigma = v\,dx + \sigma_s$ of $[M,\infty]$ if and only if  $\mf_s(\R_+)<\infty$.
\end{Prop} 
\beginpf For given $M$, we have $t_n = n$ and $M(t_{n+2}) - M(t_{n}) = 2 + \mf_s(n, n+2]$, hence 
$$
\widetilde\K[M, \infty] = \sum_{n \ge 0}(2 \cdot (2 + \mf_s(n, n+2]) - 4) = 2\sum_{n \ge 0}\mf_s(n, n+2].
$$
It remains to use Theorem \ref{t2}. \qed

\medskip

 The next result shows that logarithmic integral can converge even if $\mf_s(\mathbb{R}_+)=\infty$. 

\medskip

\begin{Prop} There exists a string $[M, L]$ with $L<\infty$ and $\mf_s[0,L) = +\infty$ such that 
$$
\int_0^\infty \frac{\log v(x)}{\sqrt x(x+1)}\,dx> -\infty
$$ 
for its spectral measure $\sigma = v\,dx + \sigma_s$.
\end{Prop} 
\beginpf Consider any sequence $\{\eps_n\} \subset(-1, 1)$, and define $\delta_{t_n}= \prod_{j=0}^n(1+\eps_j)$, $t_0=0, t_n= \sum_{j=0}^{n-1}\delta_{t_j}$ for integer $n \ge 0$, and let  $L = \sup_{n \ge 0} t_n$. Consider the function
\[
M'(t)= M_n= (\delta_{t_n})^{-2}, \qquad t\in [t_n, t_{n+1}], \qquad n \ge 0.
\]
Define the measure $\mf$ by $\mf= M'dt+\mf_s$, where $\mf_s$ is some singular measure, and let $M(t) = \mf[0,t]$ for $t\ge 0$. Then, the condition \eqref{sdk7} for $[M,L]$ is satisfied if and only if
\begin{equation}\label{sdgj}
\left\{(\delta_{t_n}+\delta_{t_{n+1}})\left(\frac{1}{\delta_{t_n}}+\frac{1}{\delta_{t_{n+1}}}\right)-4\right\}\in \ell^1
\end{equation}
and 
\begin{equation}\label{sdgh}
\left\{(\delta_{t_n}+\delta_{t_{n+1}})(\Delta \mf_s)_n\right\}\in \ell^1,
\end{equation}
where $(\Delta \mf_s)_n = \mf_s(t_n, t_{n+2}]$ for $n \ge 0$. Condition \eqref{sdgj} is satisfied if and only if
\[
\left\{(1+\eps_n)+(1+\eps_n)^{-1}-2\right\}\in \ell^1,
\]
or, equivalently, $\{\eps_n\}\in \ell^2$. If we choose $\eps_n= -(n+1)^{-\alpha},\alpha\in (\tfrac{1}{2},1)$, then
$
\sum_{n=1}^\infty (t_{n+2} - t_n)<\infty
$
and we have $L<\infty$. Condition \eqref{sdgh} in that case can be satisfied even if $\sum_n (\Delta\mf_s)_n$ diverges, that is, $\mf_s[0,L)=\infty$. For instance, we can take a singular measure $\mf_s$ such that $(\Delta \mf_s)_n = 1$ for all integers $n \ge 0$. \qed

\section{Appendix}\label{sa}
\noindent{\bf Proof of Lemma \ref{l5}.} Differentiate the function $M: r \mapsto  \idm - zJ \int_0^t \Hh(\tau)d\tau$ and use the fact that the solution to Cauchy problem \eqref{eq1} is unique. \qed

\medskip

\noindent{\bf Proof of Lemma \ref{l11}.} Put $\sigma_1 = \left(\begin{smallmatrix}0&1\\1&0\end{smallmatrix}\right)$ and $M_{\sigma_1} = \sigma_1 M \sigma_1$, where $M$ is the solution of \eqref{eq1}. Using identity $\sigma_1 \Hh \sigma_1 = J^*\Hh J =  \Hh_d$ and $J \sigma_1 = -\sigma_1 J$, it is easy to check that $JM_{\sigma_1}' = - z \Hh_d M_{\sigma_1}$. It follows that $M_{\sigma_1}(t,z) = M^d(t,-z)$ for all $t \ge 0$, $z \in \C$. Using \eqref{eq69}, we get
$$
\begin{pmatrix}\Phi^-(t,z) &  \Theta^-(t,z) \\ \Phi^+(t,z) & \Theta^+(t,z)\end{pmatrix} = \begin{pmatrix}\phantom{-}\Phi^-(t,-z) &  -\Theta^-(t,-z) \\ -\Phi^+(t,-z) & \phantom{-}\Theta^+(t,-z)\end{pmatrix}
$$
for all $t \ge 0$ and $z \in \C$. From \eqref{eq9}, one has $m(z) = -m(-z)$ for $z \in \C\setminus \R$, hence 
$$
\frac{1}{\pi}\int_{\R_+}\frac{\Im z}{|x-z|^2}\,d\mu(x) + b \Im z = \frac{1}{\pi}\int_{\R_+}\frac{\Im z}{|x+z|^2}\,d\mu(x) + b \Im z, \qquad z \in \C^+.
$$ 
This implies that $\mu$ is even.  Using $m(i+1) = -m(-i-1)$, we conclude that $a = 0$. 

\medskip

Conversely, suppose that $\mu$ is even and $a = 0$. 
The approximation procedure in Section 9 of \cite{Romanov} gives a sequence of even measures $\mu_N$ supported at finitely many points such that the corresponding Hamiltonians, $\Hh_N$, constructed in Theorem 7 of \cite{Romanov} are diagonal and $\lim_{N \to \infty}\bigl\|\int_{0}^{t}(\Hh_N(s) - \Hh(s))\,ds\bigr\| = 0$ for every $t \ge 0$. It follows that $\Hh$ is diagonal, as required. \qed

\medskip

\noindent{\bf Proof of Lemma \ref{l12}.} Let $\Hh$ be a singular nontrivial Hamiltonian on $\R_+$ such that $(0,\eps)$ is the indivisible interval of type $\pi/2$ for some $\eps>0$. Then,
for all $z \in \C^+$, we have
\begin{equation}\label{eq67}
m(z) = \frac{\Phi^+(\eps,z) + m_\eps(z)\Phi^-(\eps,z)}{\Theta^+(\eps,z) + m_\eps(z)\Theta^-(\eps,z)} 
= z \int_{0}^{\eps}\langle \Hh(t) \zo, \zo\rangle\,dt + m_\eps(z),  
\end{equation}
by formula \eqref{eq10} for $r = \eps$ and Lemma \ref{l5}. So, we have $b \ge  \int_{0}^{\eps}\langle \Hh(t) \zo, \zo\rangle\,dt$ in this situation. 

\medskip

Conversely, assume that $b> 0$ in \eqref{eq371}. Consider a Hamiltonian $\Hh_{(b)}$ whose Weyl-Titchmarsh function $m_{\Hh_{(b)}}$ coincides with $m - bz$. Define 
$$
\widetilde \Hh(x) = \begin{cases}\diag(0,1), &x \in [0, b], \\ \Hh_{(b)}(x-b), &x > b.\end{cases}
$$
Let $m_{\widetilde \Hh}$ denote the Weyl-Titchmarsh function of $\widetilde\Hh$. Then, a variant of \eqref{eq67} for $\widetilde \Hh$, $\eps = b$, gives
$$
m_{\widetilde\Hh} = bz + m_{\Hh_{(b)}} = bz + m -bz = m.
$$
Thus, the Weyl-Titchmarsh functions of $\Hh$ and $\widetilde\Hh$ coincide. It follows from de  Branges theorem formulated in the Introduction that the Hamiltonians $\Hh$, $\widetilde \Hh$ are equivalent. Hence, there is an absolutely continuous strictly increasing function $\eta \ge 0$ such that $\widetilde\Hh(t) = \eta'(t)\Hh(\eta(t))$ almost everywhere on $\R_+$. In particular, the interval $(0, \eta(b))$ is indivisible of type $\pi/2$ for $\Hh$. It follows that for $\eps = \eta(b)$ we have
$$
b = \int_{0}^{b}\trace{\widetilde\Hh(t)}\,dt = \int_{0}^{\eta(b)}\trace{\Hh(s)}\,ds = \int_{0}^{\eps}\langle\Hh(s)\zo,\zo\rangle\,ds, 
$$ 
completing the proof of the lemma. \qed

\medskip 

\noindent{\bf Proof of Lemma \ref{l9}.} The matrix-function   
$$
M(t,z) = \begin{pmatrix}\phantom{-}   \cos(t\sqrt{a_1 a_2} z) &\sqrt{a_2/a_1}\sin(t\sqrt{a_1 a_2}z)\\ -\sqrt{a_1/a_2}\sin(t\sqrt{a_1 a_2}z) & \cos(t\sqrt{a_1 a_2}z)\end{pmatrix}  
$$
solves Cauchy problem \eqref{eq1} for $\Hh = \diag(a_1,a_2)$. It follows from \eqref{eq9}  that the Weyl-Titchmarsh function of $\Hh$ is given by $m(z) = i\sqrt{a_2/a_1}$ for all $z \in \C^+$. Taking imaginary part, we get $w_r(x) = \sqrt{a_2/a_1}$, $x \in \R$, and $\log\I_{\Hh}(r) = \J_{\Hh}(r) = \log\sqrt{a_2/a_1}$ for all $r \ge 0$, as required. \qed\bigskip

\bibliographystyle{plain} 
\bibliography{bibfile}
\enddocument